\documentclass[12pt]{amsart}
\usepackage[latin1]{inputenc}
\usepackage{amssymb,pstricks}
\usepackage{graphicx}
\newtheorem{thm}{Theorem}[section] 
\newtheorem{pro}[thm]{Proposition}
\newtheorem{lem}[thm]{Lemma}

\theoremstyle{definition}
\newtheorem{dfn}{Definition}[section]

\theoremstyle{remark}
\newtheorem*{rem}{Remark}
\newtheorem*{rems}{Remarks}
\newcommand{\Z}{\mathbb{Z}}
\newcommand{\V}{\mathbf{V}\hspace{-2.5pt}}

\newcommand{\ra}{\rangle}\newcommand{\la}{\langle}

\newcommand{\fk}{\mathbf{k}}
\newcommand{\esp}{\vspace{5pt}}
\newcommand{\one}{{ \rm \setlength{\unitlength}{1em}
\begin{picture}(0.7,0.65)
\put(0,0){$1$}\put(0.34,0){\line(0,1){0.65}}
\end{picture} }\!{}}

\newcommand{\oneun}{\one_1}

\begin{document}
\title[Spin Verlinde formulas]{A spin decomposition
of the Verlinde formulas for type A modular categories}

\author{Christian Blanchet}
\address{L.M.A.M., Universit\'e de Bretagne-Sud,
BP 573, F-56017 Vannes, France  }
\email{Christian.Blanchet@univ-ubs.fr}
\date{February 2005}

\begin{abstract}
 A modular category is a braided category with some additional
algebraic features.
The interest of this concept is that it provides
a Topological Quantum Field Theory in dimension $3$.
 The Verlinde formulas associated with a modular category are the dimensions
 of the TQFT modules.
  We  discuss reductions and
  refinements of these formulas for   modular categories related with $SU(N)$.
 Our main result is a splitting of the Verlinde formula,
corresponding to a brick decomposition of the TQFT modules
whose summands are indexed by spin structures modulo
an even integer.
 We introduce here the notion of a spin modular category, and give 
the proof of the decomposition theorem in this general context.  
\end{abstract}

\maketitle\setcounter{section}{-1}
\section{Introduction}
Given a simple, simply connected complex Lie group $G$,
 the Verlinde formula \cite{Ver}
    is a combinatorial function $\mathcal{V}_G: 
(K,g)\mapsto \mathcal{V}_G(K,g)$
associated
    with $G$ (here the integers $K$ and $g$ are respectively
the level and the genus). In conformal field theory this formula
gives the dimension of the so called { conformal blocks}.
 Its combinatorics was intensively studied 
since this formula  has a deep interpretation 
as the rank of a space of { generalized theta functions}
(sections of some bundle over the moduli space of $G$-bundles
over a Riemann surface) \cite{BL,Beau,Fal,So}.
See \cite{BiL,BiL2}, for a  development using
methods of sympleptic geometry.

We will consider here  a purely topological approach
to Verlinde formulas related with  $SU(N)$.
 The genus $g$ Verlinde formula associated with a modular category 
\cite{Tu} is
 the dimension
 of the TQFT-module of a genus $g$ surface; the general formula is
given in \cite[IV,12.1.2]{Tu}.
 Various constructions of modular categories are known, either
from quantum groups \cite{AP,BK,Sa} or from skein theory \cite{TW2,Bhec,BB}.  
The geometric Verlinde formula for the group $SU(N)$
at level $K$ is recovered from the so called $SU(N,K)$
modular category. This modular category can be 
 obtained either
from the quantum group $U_qsl(N)$ when $q=s^2$ is a primitive
$(N+K)$-th root of unity or from Homfly skein theory. Its simple objects correspond to
the weights in the fundamental alcove.
 
 One may also consider a modular category with less simple objects.
This was done for $gcd(N,K)=1$ by restricting to representations whose heighest weight is in the root lattice, and was called the  projective or $PSU(N)$ theory \cite{KT,Yo,MW,Le,LT,Sa2}.
 Using an appropriate choice of the framing parameter in Homfly skein theory,
we have obtained in \cite{Bhec} a variant which is defined
for all $N,K$. We are not aware of a quantum group approach
to these reduced modular categories for $\gcd(N,K)>1$. Nevertheless we find it convenient to call them $PU(N,K)$ modular categories. In our construction
the simple object corresponding to the deformation of the determinant of
the vector representation of $sl(N)$ may be non trivial; we think
that a version of the quantum group $U_q(gl(N))$ could be used here. 

As well known, the Verlinde formula  for the $SU(N,K)$ modular category
coincides with the formula in conformal fields theory for the group $SU(N)$;
 $$d_{N,K}(g)=\mathcal{V}_{SU(N)}(K,g).$$
 We show that for the $PU(N,K)$ modular category 
 the Verlinde formula is
 $$\tilde{d}_{N,K}(g)=\frac{d_{N,K}(g)}{{N'}^g},$$
 where $N'=\frac{N}{\gcd(N,K)}$.
 These integral numbers satisfy the level-rank duality relation
 $$\tilde{d}_{N,K}(g)=\tilde{d}_{K,N}(g)\ ,$$
 which is an integral version
 of a reciprocity formula in \cite{rec} (see also \cite{KT}).

 Our main contribution here is to show that under certain condition the TQFT modules decompose in blocks indexed by spin type structures (respectively $1$-dimensional
 cohomology classes) on the surface, and compute the corresponding refined
 Verlinde formulas. 
 
An important part of this paper is devoted to the spin decomposition theorem.
 The proof  is given in the general case of a 
modulo $d$ spin  modular category; this notion, developed in Section
\ref{spin-g}, is new and appears in the  $\mathbb{Z}/d$ graded cases which are  
not weakly non-degenerate in \cite{LT}. 
As a motivation, we give below the combinatorial counterpart
of  this theorem for the A series, in the special case where the rank
is even and 
divides the level (Theorem \ref{spin_A}).
 We consider the action of $\mathbb{Z}/N$ on the set
$$\Gamma_{N,K}=\{\lambda=(\lambda_1,\dots,\lambda_{N}),\
K\geq\lambda_1\geq\dots\geq \lambda_{N-1}\geq\lambda_N=0\}\ ,$$ 
given for the generator 
 of the cyclic group
$\mathbb{Z}/N$ by
$$(\lambda_1,\dots,\lambda_{N-1},0)\longmapsto
(K,\lambda_1,\dots,\lambda_{N-1})-(\lambda_{N-1},\dots,\lambda_{N-1})\ .$$ 
 We denote by $\sharp\, \mathrm{orb}(\lambda)$
 the cardinality of the orbit of $\lambda$,
and by $Stab(\lambda)$ the stabilizer subgroup. 
 For $\mathfrak{a},\mathfrak{b}\in\mathbb{Z}/N$, the numbers 
$\epsilon_\lambda(\mathfrak{a},\mathfrak{b}) \in
\{0,1,-\frac{1}{2},\frac{1}{2}\}$ are defined  as follows.

\esp\noindent If $\sharp\, \mathrm{orb}(\lambda)$ is even, then
$$\epsilon_\lambda(\mathfrak{a},\mathfrak{b})=\left\{\begin{array}{l}
1\ \text{ if  $\mathfrak{a}$ and $\mathfrak{b}$  are zero modulo $|Stab(\lambda)|$,}\\
0 \ \text{ else.}\end{array}\right.$$
If $\sharp\, \mathrm{orb}(\lambda)$ is odd, then
$$\epsilon_\lambda(\mathfrak{a},\mathfrak{b})=\left\{\begin{array}{l}
\frac{1}{2}(-1)^{\frac{2\mathfrak{a}}{|Stab(\lambda)|}\frac{2\mathfrak{b}}{|Stab(\lambda)|}}
\ \text{ if  $\mathfrak{a}$ and $\mathfrak{b}$ 
are  zero modulo $\frac{|Stab(\lambda)|}{2}$,}\\
0\ \text{ else.}\\
\end{array}\right.$$

\begin{thm} Suppose that
$N$ is even, and that $K/N$ is an odd integer.\\
a) For $(a,b)\in(\mathbb{Z}/N)^g\times (\mathbb{Z}/N)^g$, the formula 
\begin{eqnarray*} d_{N,K}^{(a,b)}(g)&=&
\left((N+K)^{N-1}N\right)^{g-1}
  \sum\limits_{\lambda\in\Gamma_{N,K}}\ \prod\limits_{\nu=1}^g
\frac{\epsilon_\lambda(a_\nu,b_\nu)}{(\sharp\, \mathrm{orb}(\lambda))^2}\\
& &\times\prod\limits_{1\leq i<j\leq N}
  \left({2\sin{(\lambda_i-i-\lambda_j+j)\frac{\pi}{N+K}}}
\right)^{2-2g} 
\ .\end{eqnarray*}
defines a natural number $d_{N,K}^{(a,b)}(g)$,\\
b) There exists a splitting of the  $SU(N)$ Verlinde formula at level $K$
$$ \mathcal{V}_{SU(N)}(K,g)=\sum\limits_{(a,b)\in (\mathbb{Z}/N)^g\times (\mathbb{Z}/N)^g}
d_{N,K}^{(a,b)}(g)\ .$$
\end{thm}

For $N=2$, the spin TQFT producing the above decomposition 
was studied in \cite{BM}, and a nice 
algebrico-geometric interpretation  was obtained by 
Andersen and Masbaum \cite{AM}.
We quote that for $N>2$, the involved spin structures are not the usual ones.
 These structures  have coefficients modulo an even integer; they
can be understood as something intermediate between usual spin structures
(with modulo $2$ coefficients) and complex spin structures.
 The convenient formalism for the TQFT involving these structures
should be a slightly extended version of Homotopy Quantum Field Theory
as developed by Turaev \cite{Tu3,Tu4}.

The paper is organized as follows. In Section 1 we study spin structures
modulo an even integer. In Section 2 we define our spin modular categories.
 In Section 3 we establish the spin decomposition of the TQFT in a general
context. In Section 4 we consider Verlinde formulas for modular
categories of the A series. 
In Section 5 we establish similar decomposition theorems based
on $1$-dimensional cohomology classes. In Section 6 we give computer results for small
values of $N$ and $K$. 
\section{Spin structures modulo an even integer}
Let $d$ be an even integer. We recall 
 here the topological definition for modulo $d$  spin  structures
that was given in \cite{varso,Bhec}.

There exists, up to homotopy, a
unique non-trivial map
$g$ from the classifying space $BSO$ to the Eilenberg-MacLane space $K(\mathbb{Z}/d,2)$.
 Define the fibration
$$\pi_d: BSpin(\mathbb{Z}/d)\rightarrow BSO$$
to be the pull-back, using $g$, of the path
fibration over $K(\mathbb{Z}/d,2)$.
 The space $BSpin(\mathbb{Z}/d)$ is a classifying space
for the non-trivial central extension
of the Lie group $SO$ by $\mathbb{Z}/d$,
which we denote by $Spin(\mathbb{Z}/d)$. For $d=2$, this group
$Spin(\mathbb{Z}/2)=Spin$ is the universal cover of $SO$, and for general
$d$, we have
$$ Spin(\mathbb{Z}/d)=\frac{Spin\times \mathbb{Z}/d}{(-1,d/2)}\ .$$

Now we can use the fibration $\pi_d$ to define structures.  Let \mbox{$E_{Spin(\mathbb{Z}/d)}=\pi_d^*(E_{SO})$}
be the pull-back of the canonical vector bundle over $BSO$.

\begin{dfn} A   modulo $d$ spin structure 
(or $Spin(\Z/d)$  structure) on a manifold $M$
is an homotopy class
of fiber maps from the stable tangent bundle $\tau_M$ to
$E_{Spin(\mathbb{Z}/d)}$.
\end{dfn}
If non-empty the set of these structures, denoted by
$Spin(M;\mathbb{Z}/d)$, is affinely isomorphic to
$H^1(M;\mathbb{Z}/d)$, by obstruction theory.
 Moreover the obstruction for existence is a class
 $w_2(M;\mathbb{Z}/d)\in H^2(M;\mathbb{Z}/d)$, which is the image of
 the Stiefel-Whitney class $w_2(M)$ under the homomorphism
 induced by the inclusion of coefficients
 $\mathbb{Z}/2\hookrightarrow \mathbb{Z}/d$.
 The Stiefel-Whitney class $w_2(M)$ is zero for every compact oriented
manifold whose dimension is lower or equal to $3$, hence
 spin  structures modulo $d$ exist on these manifolds.

The various descriptions of  usual spin structures \cite{Milnor} apply
to modulo $d$ spin structures.
 The above definition defines, up to equivalence, a $Spin(\mathbb{Z}/d)$ principal
bundle over the stable oriented framed bundle $PTM$ (with fiber $SO$)
whose restriction to the fiber is equivalent to
the cover map $Spin(\mathbb{Z}/d)\rightarrow SO$.
 The cover of $PTM$ defined by the modulo $d$ spin structure
is classified by a cohomology class $\sigma\in H^1(PTM,\mathbb{Z}/d)$
whose restriction to the fiber is non-trivial. The above correspondence
is   one to one;
this gives an
alternative definition, and we will identify
$Spin(M;\mathbb{Z}/d)$ with the corresponding affine sub-space of
$H^1(PTM,\mathbb{Z}/d)$.
\begin{dfn}[Alternative definition of modulo $d$ spin structures]
A   modulo $d$ spin structure 
  on an oriented manifold $M$
is  a cohomology class $\sigma\in H^1(PTM,\mathbb{Z}/d)$
whose restriction to the fiber is non-trivial.
\end{dfn}
 Observe that a spin structure can be evaluated on
a framed $1$-cycle in the manifold.

Let us consider an oriented surface $\Sigma$.
 An immersed curve has a preferred framing defined by using the tangent
vector. If a closed embedded curve $\gamma$ bounds a disc, then the evaluation
of a modulo $d$ spin structure on the corresponding framed
$1$-cycle $\tilde{\gamma}$ is $\frac{d}{2}$.
 Following \cite{A,J}, we get the theorem below which gives a convenient
description of modulo $d$ spin structures on the oriented surface $\Sigma$.
\begin{thm}
a) Let $\gamma$ denotes an embedded closed curve with $\sharp\gamma$ components.
 The assignement $\gamma\mapsto\sigma(\tilde{\gamma})+
(\sharp\gamma)\frac{d}{2}$
extends to a well defined map 
$q_\sigma: H_1(\Sigma,\mathbb{Z}/d)\rightarrow \mathbb{Z}/d$.\\
b) The map  $\sigma\mapsto q_\sigma$ defines
a canonical bijection between $Spin(\Sigma,\mathbb{Z}/d)$
and the set of maps
$q:H_1(\Sigma,\mathbb{Z}/d)\rightarrow \mathbb{Z}/d$
 such that for all $x$, $y$ one has
\begin{equation} \label{int}q(x+y)=q(x)+q(y)+\frac{d}{2}x.y\ .\end{equation}
Here $.$ denotes the intersection form
on $H_1(\Sigma,\mathbb{Z}/d)$.
\end{thm}
\begin{proof}
 The formula $\sigma(\tilde{\gamma})+(\sharp\gamma)\frac{d}{2}$ 
is unchanged if we add or remove to the embedded curve $\gamma$
a trivial component.
 Let us denote by $\gamma$ (resp. $\gamma'$) the left handed (resp. right handed) curve
 in the band move represented in Figure \ref{bandmove}.
We have that $\sharp\gamma'-\sharp\gamma=\pm 1$. By considering the Gauss map,
we see that the  cycle $\tilde \gamma' -\tilde \gamma$ is homologous
in $PT\Sigma$
to $\tilde u$ where $u$ bounds a disc  on the surface.
 We get that
 the formula is also
unchanged under this band move.
 We deduce that homologous curves in $\Sigma$ give the same result; hence we have that
$q_\sigma$ is well defined on $H_1(\Sigma,\mathbb{Z})$.
 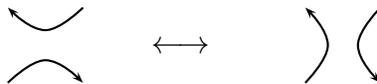
\begin{figure}
\centerline{
\begin{pspicture}[.6](0,-0.5)(1,1)
\pscurve{->}(0,0)(.5,.3)(1,0)
\pscurve{<-}(0,1)(.5,.7)(1,1)
\end{pspicture}\hspace{18pt}
$\ \longleftrightarrow\ $
\begin{pspicture}[.6](0,-0.5)(2,1)
\pscurve{->}(1,0)(1.3,.5)(1,1)
\pscurve{<-}(2,0)(1.7,.5)(2,1)
\end{pspicture}}
\caption{\label{bandmove}Band move}
\end{figure}
 
Let $\gamma$ be  a generic immersed curve. 
 Smoothing a crossing changes $\sharp\gamma$
by $\pm 1$ and does not change the $1$-cycle $\tilde \gamma$. Hence one has $q_\sigma([\gamma])=
\sigma(\tilde{\gamma})+(\sharp\gamma+I(\gamma))\frac{d}{2}$,
where $I(\gamma)$ is the number of 
double points. 
 It follows that for any $x,y\in H_1(\Sigma,\mathbb{Z})$, 
Property  (\ref{int}) holds.
 We deduce that $q_\sigma$ is well defined on $H_1(\Sigma,\mathbb{Z}/d)$.
Bijectivity is established by using that
 the map  $q_\sigma$ commutes with the action of $H^1(\Sigma,\mathbb{Z}/d)$.
\end{proof}

 Let $M=\mathbf{S}^3(L)$ be obtained by surgery on the framed link
$L$ in the $3$-sphere.
 We want to give a combinatorial description
for modulo $d$ spin structures on $M$.  Recall that $M$ is the boundary of a
$4$-manifold $W_L$ called the trace of the surgery.
To each $\sigma\in Spin(M;{\mathbb{Z}/d})$ is associated
a relative obstruction $w_2(\sigma;{\mathbb{Z}/d})$
in $H^2(W_L,M;\mathbb{Z}/d)$. 
 The group $H^2(W_L,M;\mathbb{Z}/d)$ is a free $\mathbb{Z}/d$ module
of rank $m=\sharp L$.
 Taking the coordinates of the relative obstruction
 in the preferred basis (the basis which is Poincaré dual to the cores of the handles), 
we get a map $\psi_L: Spin(M;{\mathbb{Z}/d})
\rightarrow\left(\mathbb{Z}/d\right)^m$.

The following theorem is proved in \cite{Bhec}.{Here  $B_L=(b_{ij})$ is
the linking matrix.}
\begin{thm} The map $\psi_L: Spin(M;{\mathbb{Z}/d})
\rightarrow\left(\mathbb{Z}/d\right)^m$ is injective, and its image
is the set of those $(c_1,\dots,c_m)$ which are
solutions of the following $\mathbb{Z}/d$-characteristic
equation
$$B_L\left(\begin{array}{c} c_1\\ \vdots\\ c_m\end{array}\right)=
{\frac{d}{2}}\left(\begin{array}{c} b_{11}\\ \vdots\\
 b_{mm}\end{array}\right)\ \ ({\rm mod}
\ d)\ .$$
\end{thm}

\section{Spin modular categories}
 A ribbon category is a category equipped with  tensor product,
braiding, twist and duality satisfying compatibility conditions
\cite{Tu}.
 If we are given a ribbon category $\mathcal{C}$, then we can define
an invariant of links whose components are colored with  objects of
$\mathcal{C}$. This invariant extends to a representation of the
$\mathcal{C}$-colored tangle category and more generally
to a representation of the category of $\mathcal{C}$-colored ribbon graphs
\cite[I.2.5]{Tu}.
 In a ribbon category there is a notion of trace of morphisms
and dimension of objects. 
 The trace of a morphism $f$ is denoted by $\langle f\rangle$.

\centerline{\raisebox{12mm}{$\langle f\rangle=$} \begin{picture}(10,28)
\put(6,12 ){$f$}
\put(1,2){\includegraphics[height=24mm]{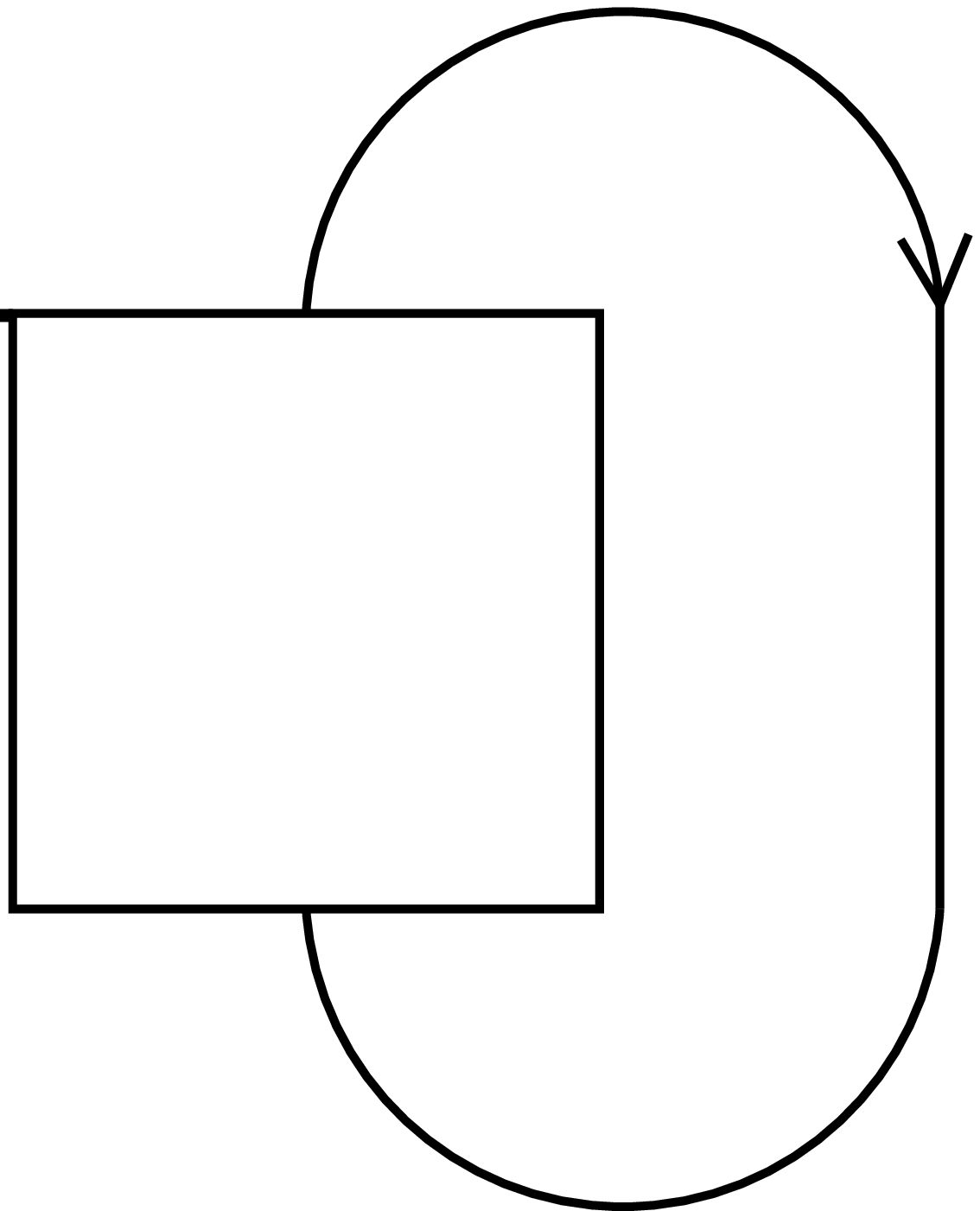}}
\end{picture}}
 The dimension
of an object $V$ is the trace of the identity morphism $\one_V$; we will use the notation
$\langle V\rangle$ as well as  $\langle \one_V\rangle$.
 We often say {\em quantum} trace and dimension
to distinguish from the usual trace and dimension in vector spaces. 

Let $\fk$ be a field. A ribbon category  is said to be
$\fk$-additive if the Hom sets are 
$\fk$-vector spaces, composition and tensor product are
bilinear, and $End(\text{trivial object})=\fk$.

We first recall the definition of a modular category
\cite{Tu,PS}. 
A {\em modular category} 
over $\fk$ 
is a $\fk$-additive ribbon category in which there
exists a finite family $\Gamma$ of simple objects $\lambda$
(here simple means that $u\mapsto u \one\!_{\lambda}$
from $\fk=End(\mathrm{trivial\ object})$ to $End(\lambda)$ is an isomorphism)
satisfying the  axioms below. 
\begin{itemize}
\item (Domination axiom) For any object $V$ in the category
there exists a finite decomposition
$\one\!_V=\sum_i f_i\one\!_{\lambda_i}g_i$,
with $\lambda_i\in\Gamma$ for every $i$.
\item (Non-degeneracy axiom) The following matrix is invertible.
$$S=(S_{\lambda\mu})_{\lambda,\mu\in \Gamma}$$
where $S_{\lambda\mu}\in \fk$ is the endomorphism
of the trivial object
associated with the $(\lambda,\mu)$-colored,
$0$-framed Hopf link  with linking $+1$. 
\end{itemize}
It follows that $\Gamma$ is a representative
set of isomorphism classes of simple objects; note that the trivial
object $\Theta$ is simple, so that we may suppose that $\Theta$
is in $\Gamma$.
 If we replace the non-degeneracy axiom
by the non-singularity condidition below
then  we have the definition of a 
{\em pre-modular category} (
a morphism $f\in Hom(V,W)$ is called negligible if
for any $g\in Hom(W,V)$ we have $\la fg\ra=0$): 
\begin{itemize}
\item  (Non-singularity) The category has no non-trivial negligible morphism.
\end{itemize}
 
 A general modularization procedure for pre-modular categories, 
and a criterion for existence are developed by Brugui\`eres 
\cite{Bru}, and by M\"uger in the context of $*$-categories
\cite{Mu}.
 Note that after quotienting by negligible we get the
non-singularity condition. This property gives that the 
pairing 
$$ \begin{array}{ccc}
Hom(V,W)\otimes  Hom(W,V)&\mapsto &\fk\\
f\otimes g&\mapsto&\la fg\ra
\end{array}$$
is non singular. We can deduce that there exists no non-trivial morphism
between non-isomorphic simple
objects.

  One may ask further that the category has direct sums.
 In fact direct sums may be added in  a formal way,
and a pre-modular category with direct sums is abelian. This latter
fact was pointed out to us by  Bruguières.

 In a modular category $\mathcal{C}$, with representative set of simple objects
$\Gamma$, the {\em Kirby color} $\Omega=\sum_{\lambda\in\Gamma}\langle \lambda \rangle \lambda$
 is used to define an invariant of  closed
oriented manifolds  with colored graph.
 If $M=S^3(L)$ is obtained by surgery on the framed link
$L$ in the sphere and contains a colored graph $K$, then a formula for this invariant is
$$\tau_{\mathcal{C}}(M,K)=
\frac{{\langle L(\Omega,
\dots,\Omega),K\rangle}}{\langle U_1(\Omega)\rangle^{b_+}
\langle U_{-1}(\Omega)\rangle^{b_-}}\ .
$$
Here $b_+$ (resp. $b_-$) is the number of
positive (resp. negative) eigenvalues of the linking
matrix $B_L$, and $U_{\pm1}$ denotes the unknot with framing $\pm1$.

Modular $G$-categories, with $G$ a group have been introduced
by Turaev in \cite{Tu4}; details in the case of an abelian group $G$,
 and examples derived
from quantum groups are given in \cite{LT}.

Let $G$ be an abelian  group. A $G$ {\em grading} of a $\fk$-additive
 monoidal category $\mathcal{C}$
is a family of full sub-categories $\mathcal{C}_j,\ j\in G$,
such that \begin{itemize}
\item[(i)] for any pair of objects $V\in Obj(\mathcal{C}_j)$,
$V'\in Obj(\mathcal{C}_{j'})$, one has  $V\otimes V'\in Obj(\mathcal{C}_{j+j'})$;\\
\item[(ii)] if for some pair of objects $V\in Obj(\mathcal{C}_j)$,
$V'\in Obj(\mathcal{C}_{j'})$ one has   $Hom_{\mathcal{C}}(V,V')\neq\{0\}$,
then $ j=j'$;\\
\item[(iii)] each object of $\mathcal{C}$ is either in  $\cup_j Ob(\mathcal{C}_j)$,
or  a direct sum of objects in $\cup_j Obj(\mathcal{C}_j)$.
\end{itemize}

Axiom (iii) asks that every object splits as a direct sum of homogeneous
objects. Axiom (i) asks that tensor product is homogeneous, and axiom (ii) that
any non-zero morphism with source or target
an homogeneous object is homogeneous; this implies that the dual of an homogeneous object 
has opposite grading.

\begin{figure}
\centerline{
\begin{picture}(60,40)
\put(0,0){\includegraphics[height=40mm]{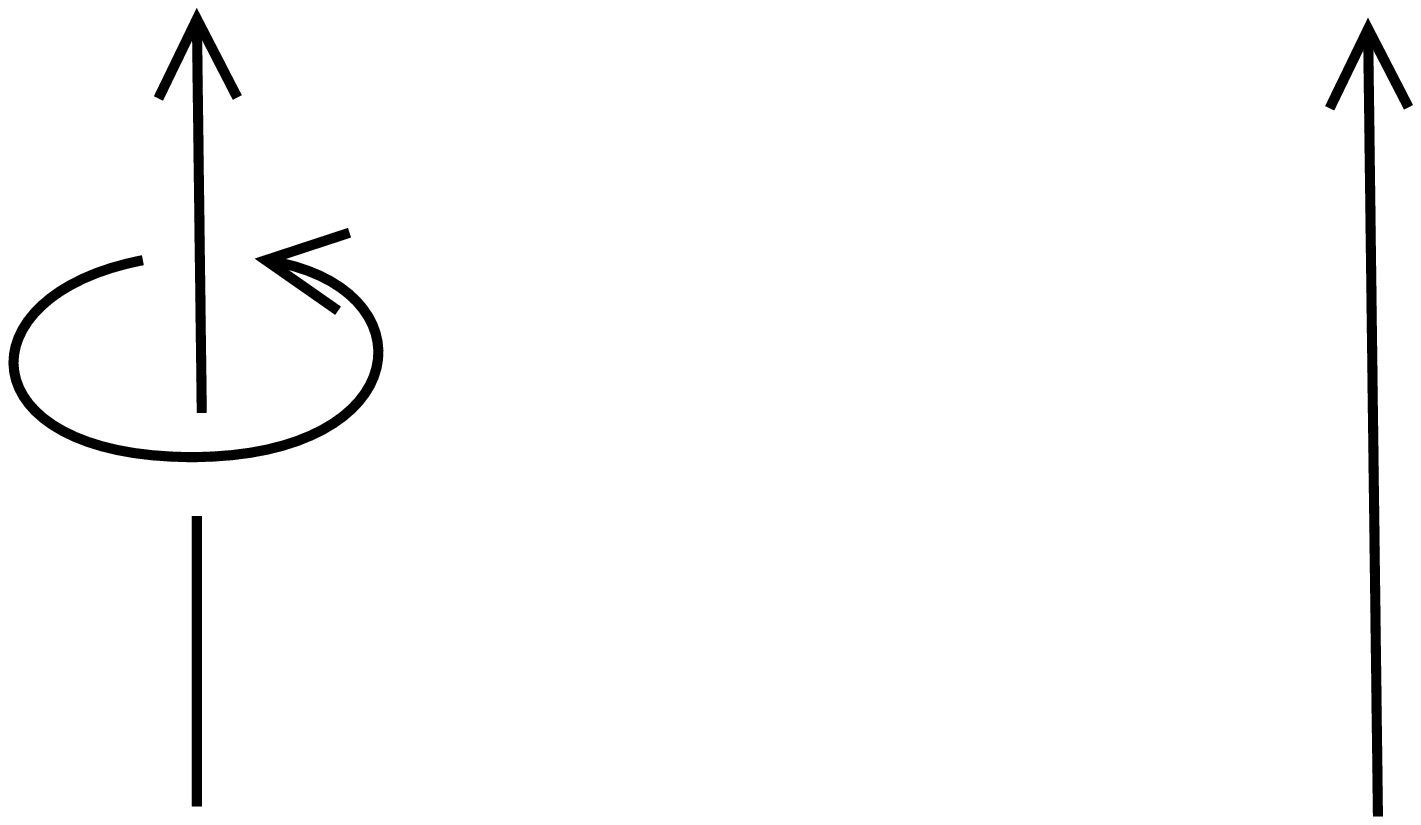}}
\put(9,5){$\lambda$}
\put(21,18){$J$}
\put(35,20){$=\ \ \ \chi(J)$}
\put(57,5){$\lambda$}
\end{picture}
}
\caption{\label{grading}}
\end{figure}

Let $\mathcal{C}$ be a modular category. We denote
by $\mathcal{U}(\mathcal{C})$ the abelian group of isomorphism
classes of invertible objects in $\mathcal{C}$
(the law is tensor product).
 If $U$ is a subgroup of ${\mathcal{U}}(\mathcal{C})$ and 
$G=\hat U$ is
the group of characters $\chi:U\rightarrow \fk^*$, then  the category
$\mathcal{C}$ is  $G$ graded. A simple object $\lambda$
is an object in $\mathcal{C}_\chi$
 if and only if for every $J\in U$
equality 
in Figure \ref{grading} holds.

A {\em modular $G$-category} \cite{LT}
over $\fk$ 
is a $G$ graded  $\fk$-additive ribbon category 
$(\mathcal{C};\mathcal{C}_j,\ j\in G)$
 in which there
exists  finite families $\Gamma_j\subset Ob(\mathcal{C}_j)$, $j\in G$,
 of simple objects $\lambda$
satisfying the  axioms below.
\begin{itemize}
\item (Domination axiom) For any object $V$ in $\mathcal{C}_j$
there exists a finite decomposition
$\one\!_V=\sum_i f_i\one\!_{\lambda_i}g_i$,
with $\lambda_i\in\Gamma_j$ for every $i$.
\item  (Non-degeneracy axiom) The following matrix is invertible.
$$S=(S_{\lambda\mu})_{\lambda,\mu\in \Gamma_0}$$
where $S_{\lambda\mu}\in \fk$ is the endomorphism
of the trivial object
associated with the $(\lambda,\mu)$-colored,
$0$-framed Hopf link  with linking $+1$. 
\end{itemize}

 It is shown in \cite{Tu4} that a modular $G$-category
with $G$ an abelian group
gives  invariants of $3$-manifolds
equipped with a $1$-dimensional cohomology class.
 A modular $G$-category may not be a modular category, even
in the case where $G$ is a finite abelian group
(see \cite[Section 1.6]{LT}).
 We point out that a modular category with a $G$ grading
is not necessarily a modular $G$-category.
 The reason is that the $S$-matrix restricted to zero graded
objects may be non-invertible. 
 In addition, the  zero graded
subcategory may be non-modularizable, so that there is no hope to get
a modular $G$-category by using some modularization procedure.
 The latter fact implies that the modular $G$-category
is not weakly non-degenerate \cite{LT};
it is verified for the class of modular categories
we consider below. These categories  have a $\Z/d$
grading with $d$ even and  give invariants of $3$-manifolds
equipped with modulo $d$ spin structures; the relevant version 
of Homotopy Quantum Field Theories as considered by Turaev,
 should be understood in relation with
\cite[Remark 7.4.6]{Tu4}. 
 
For a simple object $\lambda$ the twist coefficient $\theta_\lambda$
is defined by the figure \ref{twist}.
 In the quantum group context, this coefficient is given by the action of the 
so called quantum Casimir.
\begin{figure}
$$\begin{pspicture}[.5](0,-0.5)(1,1)
\pscurve{-}(.1,0)(.15,.45)(.45,.75)(.7,.5)(.4,.2)(.20,.38)
\pscurve{->}(.14,.62)(.1,.9)(.1,1)
 \put(-.2,0){$\lambda$}
\end{pspicture}
\ =\ \theta_\lambda\ 
\begin{pspicture}[.5](0,-.5)(1,1)
\psline{->}(.1,0)(.1,1)
 \put(.2,0){$\lambda$}
\end{pspicture}$$
\caption{\label{twist}}
\end{figure}
\label{spin-g}
\begin{dfn}
Let $d$ be an even integer (resp. an integer). 
A  modular category is modulo $d$ spin\footnote{Invariants associated with a (mod. 2) spin modular category
are considered in \cite{Sa3}[Theorem 2b]. Our definition here is more general.} (resp. modulo $d$ cohomological) 
if it is equipped
with
 an invertible object $\varrho$ whose order is $d$ and whose twist coefficient
is $\theta_\varrho=-1$ (resp. $\theta_\varrho=1$).
\end{dfn}

In the following we will mainly discuss the spin case; the cohomological case
will be considered in section \ref{coho_}.
  
Let $(\mathcal{C},\varrho)$ be a modulo $d$ spin  modular category,
with $\Gamma$ as a representative set of simple objects. 
 The object $\varrho^d$ is isomorphic to the trivial, hence we have
$\la \varrho^{d}\ra=1$. The dual objects $\varrho$ and 
$\varrho^{d-1}$ have the same quantum dimension.
We deduce that $\mathfrak{d}=\la \varrho\ra=\pm 1$.
 Note that invertible objects are simple, hence the braiding for
$\varrho^2$ is identity up to a scalar.
 By closing we get this scalar and  establish the  following identity.\\[10pt]
\begin{equation}\label{cross_rho}
\begin{pspicture}[.5](0,-0.5)(1,1)
\psline{->}(1,0)(0,1)
\psline{-}(0,0)(.35,.35)
\psline{->}(.65,.65)(1,1)
 \put(0,.3){$\varrho$} \put(.8,.3){$\varrho$}
\end{pspicture}
\ =\ -\ 
\mathfrak{d}\ \ 
\begin{pspicture}[.5](0,-0.5)(1,1)
\pscurve{->}(0,0)(.3,.5)(0,1)
\pscurve{->}(1,0)(.7,.5)(1,1)
 \put(0,.4){$\varrho$} \put(.8,.4){$\varrho$}
\end{pspicture}
\end{equation}

The next identity is obtained in a similar way.\\[10pt]

\begin{equation}\label{band_rho}
\begin{pspicture}[.5](0,-0.5)(1,1)
\pscurve{->}(0,0)(.3,.5)(0,1)
\pscurve{<-}(1,0)(.7,.5)(1,1)
 \put(0,.4){$\varrho$} \put(.8,.4){$\varrho$}
\end{pspicture}
\ =\ \mathfrak{d}\ \ \ 
\begin{pspicture}[.5](0,-0.5)(1,1)
\pscurve{->}(0,0)(.5,.3)(1,0)
\pscurve{<-}(0,1)(.5,.7)(1,1)
 \put(.3,0){$\varrho$} \put(.4,.9){$\varrho$}
\end{pspicture}
\end{equation}
 
 It is convenient to fix a primitive $d$-th root of unity $\zeta$,
and to identify the group of characters 
$\chi :\{\varrho^j,j\in\mathbb{Z}/d\} \rightarrow \fk^*$ with $\Z/d$.
 Then the category $\mathcal{C}$ is $\Z/d$ graded.
 A simple object $V$ has degree equal to $j$ if and only if the equality
in figure \ref{varrho} holds.
\begin{figure} 
\centerline{
        \begin{pspicture}(0,-.5)(5,2)
        \put(2.1,0){$\varrho$}\put(3.0,0){$V$}
        \psline[linearc=.3](2,0)(2,1.05)
        \psline[linearc=.3]{->}(2,1.35)(2,1.7)
        \pscurve(2.8,0)(2.7,.5)(2.2,.7)
        \pscurve{->}(1.9,.8)(1.7,1)(2.75,1.4)(2.8,1.7)
\put(4,.7){ $=\ \ \zeta^{j}$}
\end{pspicture}
        \begin{pspicture}(1,-0.5)(5,2)
        \put(2.1,0){$\varrho$}\put(3,0){$V$}
        \psline[linearc=.3]{->}(2,0)(2,1.7)
        \psline{->}(2.8,0)(2.8,1.7)
        \end{pspicture}}
\caption{\label{varrho}}
\end{figure}

 The Kirby color  decomposes according to
this grading.
$$\Omega=\sum_{\lambda\in \Gamma}
\la \lambda \ra \lambda\ =\sum_{j\in \Z/d}\ \Omega_j$$
 Here the  notation $\la \lambda \ra$ is the quantum dimension of $\lambda$.

The proof of the  theorem below
is the same as in the ungraded case (see e.g. \cite{BB}). 
 The statement holds for
any  $G$ graded pre-modular category with $G$ an abelian group
\cite[Prop. 1.4]{LT}.

\begin{thm}(Graded sliding property)
Suppose that $V_j$ is an object in $\mathcal{C}_j$, then the equality
in Figure \ref{graded_sliding} holds for any $j'\in\mathbb{Z}/d$.
 Here the framed knot labeled with $\Omega_{j'}$ may be knotted or 
linked with the other component labeled $V_{j}$; this fact is represented by the dashed part in the figure.
\begin{figure} 
\centerline{\includegraphics[height=50mm]{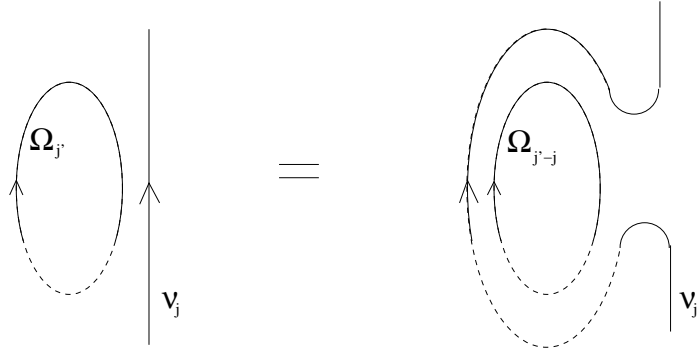}}
\caption{\label{graded_sliding}Graded sliding property.}
\end{figure}
\end{thm}

\label{spinref}
The following theorem is proved from the graded sliding property
as was done in \cite[Theorem 4.2]{Bhec}.
 We suppose 
that $\mathcal{C}$ is a modulo $d$ spin  modular
category, and that 
$$\Omega=\sum_{\lambda\in \Gamma}
\la \lambda \ra \lambda\ =\sum_{j\in \Z/d}\ \Omega_j$$
is the graded decomposition of the Kirby color;
note \cite[lemma 4.5]{Bhec} that $\langle U_{\pm 1}(\Omega)\rangle=
\langle U_{\pm 1}(\Omega_{d/2})\rangle$.
\begin{thm}Let $\mathcal{C}$ be a modulo $d$ spin modular category,
and 
\mbox{$\Omega=\sum_{j\in \Z/d}\ \Omega_j$} be the graded decomposition of the Kirby
element.
Provided $c=(c_1,\dots,c_m)$ satisfies the modulo $d$
characteristic condition,
the formula
$$\tau_{\mathcal{C}}^{\mathrm{spin}}(M,\sigma)=
\frac{{\langle L(\Omega_{c_1},
\dots,\Omega_{c_m})\rangle}}{\langle U_1(\Omega)\rangle^{b_+}
\langle U_{-1}(\Omega)\rangle^{b_-}}
$$
defines an invariant of the surgered manifold $M={\bf S}^3(L)$
equipped with the modulo $d$ spin  structure
$\sigma=\psi_L^{-1}(c_1,\dots,c_m)$
.\\
Moreover,
$$\forall M\ \ \tau_{\mathcal{C}}(M)=\sum_{\sigma\in Spin(M;\mathbb{Z}/d)}
\tau_{\mathcal{C}}^{\mathrm{spin}}(M,\sigma)$$
\end{thm}
\section{The spin decomposition of the Verlinde formula}\label{decomp}
If we are given a modular category $\mathcal{C}$ then we get
a TQFT. In brief we have a functor $\V_\mathcal{C}$
from a cobordism category in dimension $3$
to vector spaces. If $\mathcal{C}$ is a 
modulo $d$ spin modular category, then we will construct here a  decomposition
of the TQFT modules $\V_\mathcal{C}(\Sigma_g)$
of a genus $g$ surface
and compute the ranks of the summands.

The TQFT gives a normalized invariant for a closed $3$-manifold $M$
equipped with $p_1$-structure or $2$-framing $\alpha$ and colored graph $K$. 
We extend the scalar field $\fk$ if necessary, and fix $\kappa$ such that
$\kappa^6=\frac{\la U_1(\Omega)\ra }{\la U_{-1}(\Omega)\ra}$.
 Let $\mathcal{D}=\kappa^{-3}\la U_1(\Omega)\ra=
\kappa^{3}\la U_{-1}(\Omega)\ra$; note that $\mathcal{D}^2=\la \Omega \ra$.

The normalized invariant of a connected closed $3$-manifold
\mbox{$M=(M,\alpha,K)$} is then \cite{BHMV}
\begin{equation} \label{normalized} Z_\mathcal{C}(M,\alpha,K)=\mathcal{D}^{-1-b_1(M)} \kappa^{\sigma(\alpha)}
\tau_{\mathcal{C}}(M,K) \ .\end{equation}
Here $b_1(M)$ is the first Betti number, and $\sigma(\alpha)$
is the sigma invariant:
$\sigma(\alpha)=3signature(W_L)-\langle p_1(W_L,\alpha),[W_L]\rangle$,
where $W_L$ is the trace of the surgery and $p_1(W_L,\alpha)\in H^4(W_L,S^3(L))$ is the relative
obstruction to extending $\alpha$.

Let $\Sigma$ be an oriented surface with
structure (a marking \cite{Tu} or a $p_1$-structure
\cite{BHMV}). We  use the object $\varrho$
to define  a group action on $\V_\mathcal{C}(\Sigma)$ as follows.
 To an embedded oriented curve $\gamma$ in $\Sigma$ we  associate
the TQFT operator
$$\phi_\gamma :\V_\mathcal{C}(\Sigma)\rightarrow \V_\mathcal{C}(\Sigma)$$
corresponding to a trivial cobordism $[0,1]\times \Sigma$ equipped with
a colored link $\gamma_{\frac{1}{2}}(\varrho)$.
 Here $\gamma_{\frac{1}{2}}$ is the link  ${\frac{1}{2}}\times\gamma$
 equipped with the framing given by the orientation and the
normal vector parallel to $\Sigma$.
The components of this link are colored with $\varrho$.

 The spectral projector of $\phi_\gamma$ corresponding to the
eigenvalue $\zeta^\nu$ is 
equal to $\frac{1}{d}\sum_{j=0}^{d-1} \zeta^{-\nu j}\phi_\gamma^j$.
 This projector is represented by a trivial cobordism
with colored link $\gamma(\pi_\nu)$ where
the color $\pi_\nu$ is defined by
$$\pi_\nu=\frac{1}{d}\sum_{j=0}^{d-1} \zeta^{-\nu j}\varrho^j\ .$$
 Using the definition of the grading, we get the following
lemma.

\begin{figure}
\centerline{
\begin{picture}(70,45)
\put(8,7){$V$}
\put(16,20){$\pi_\nu$}
\put(35,20){$=\ \ \delta_{\nu j}$}
\put(57,7){$V$}
\put(0,5){\includegraphics[height=40mm]{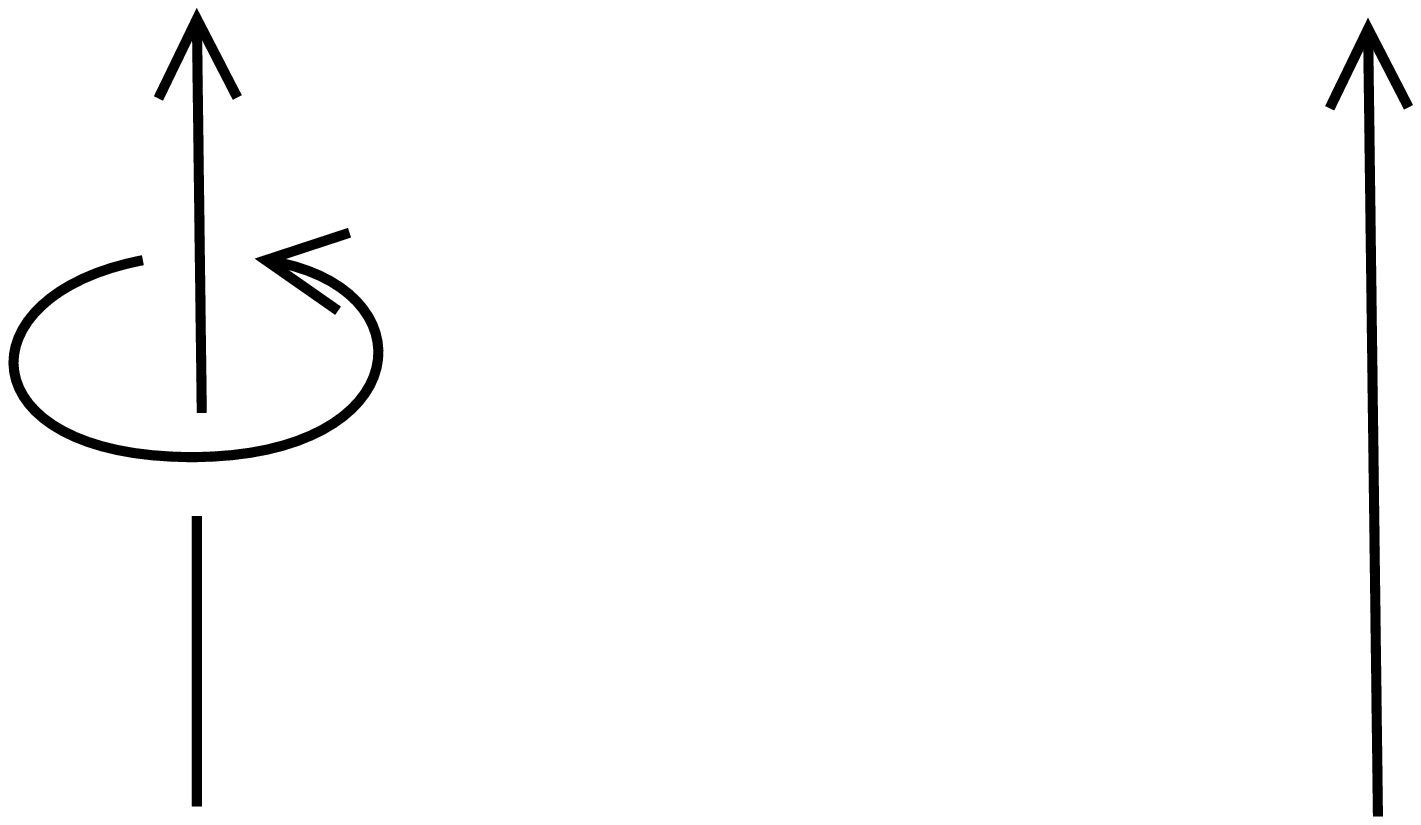}}
\end{picture}}
\caption{\label{proj_fig}}
\end{figure}

\begin{lem}\label{proj}
Let $V$ be an object in $\mathcal{C}_j$; denote by
 $\delta_{\nu j}$   the Kronecker symbol. One has the equality
in Figure \ref{proj_fig}.
\end{lem}

We denote by $PT\Sigma$ the principal $SO$-bundle
of oriented orthonormal frames in the stabilized tangent
bundle to $\Sigma$ (we could stabilize only once).

We denote by $\tilde\gamma$ the lift in $PT\Sigma$, using the unit
tangent vector,
of the embedded  curve $\gamma$.
\begin{pro}
There exists a well defined action
of the group $H_1(PT\Sigma,\mathbb{Z}/d)$ on $\V_\mathcal{C}(\Sigma)$, 
which maps $x=[\tilde\gamma]$ to the operator
$\psi_x=(-\mathfrak{d})^{\sharp\gamma}\phi_\gamma$.
\end{pro}
\begin{proof}
The $\mathbb{Z}/d$ module $H_1(PT\Sigma,\mathbb{Z}/d)$
is generated by the $1$-cycles $\tilde \gamma$ associated with embedded curves $\gamma$.
 A trivial circle represents the generator on the fiber;
 this generator has order $2$. A disjoint union represents the sum.
 All the other 
 relations are given by the modified band move
 in Figure \ref{mbandmove}.
 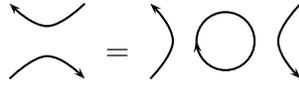
\begin{figure}
\centerline{
\begin{pspicture}[.5](0,-0.5)(1,1)
\pscurve{->}(0,0)(.5,.3)(1,0)
\pscurve{<-}(0,1)(.5,.7)(1,1)
\end{pspicture}
$\ =\ $
\begin{pspicture}[.5](0,-0.5)(2,1)
\pscurve{->}(0,0)(.3,.5)(0,1)
\pscircle(1,.5){.4}
\pscurve{->}(.63,.42)(.615,.47)(.61,.50)
\pscurve{<-}(2,0)(1.7,.5)(2,1)
\end{pspicture}}
\caption{\label{mbandmove}Modified band move}
\end{figure}
By relation (\ref{band_rho}) the modified band move doesn't change the number of components
modulo $2$, hence by relation (\ref{cross_rho}) $\psi_x$ is well defined
by the formula $\psi_x=(-\mathfrak{d})^{\sharp\gamma}\phi_\gamma$.
Here $\gamma$ is an embedded curve such that
the lift $\tilde\gamma$ represents $x$.
 A crossing resolution changes by $\pm 1$ the number of components,
 hence in the above formula we can use an immersed
 curve as well. If $\gamma,\gamma'$ represent
$x$ and $x'$, then we can isotope $\gamma'$ so that
$\gamma\cup \gamma'$ is an immersed
curve.  This shows that for all $x,x'$, one has
 $\psi_x\psi_{x'}=\psi_{x+x'}$.
 \end{proof}
As a consequence, we have a decomposition of $\V_\mathcal{C}(\Sigma)$ indexed by
the group $H^1(PT\Sigma,\mathbb{Z}/d)$ identified with the characters
on $H_1(PT\Sigma,\mathbb{Z}/d)$.
 Recall that we have chosen a primitive $d$-th root of unity
denoted by $\zeta$. A vector $v$ belongs to the component
indexed by $\sigma$ if and only if
for every $x\in H_1(PT\Sigma,\mathbb{Z}/d)$ one has
$\psi_x v=\zeta^{\sigma(x)}v$ .
 Since the generator of the fiber acts by $-1$, only the classes
whose restriction to the fiber is non-trivial,
i.e. $Spin(\Z/d)$ structures,
will correspond to non-trivial summands.

If $\sigma$ is a modulo $d$ spin structure  on the genus $g$ oriented surface
$\Sigma_g$, we denote by ${\V}_\mathcal{C}(\Sigma_g,\sigma)$
the corresponding summand and by $ d_\mathcal{C}(g,\sigma)$ its dimension.
$${\V}_\mathcal{C}(\Sigma_g,\sigma)=\{v\in {\V}_\mathcal{C}(\Sigma_g),
\forall x\in H_1(PT\Sigma,\mathbb{Z}/d)\  \psi_x v=\zeta^{\sigma(x)}v\}$$
\begin{thm}\label{spin}
a) There exists a splitting of the TQFT module
$$ \V_\mathcal{C}(\Sigma_g)=
\oplus_{\sigma\in Spin(\Sigma_g,\mathbb{Z}/d)}\ 
\V_\mathcal{C}(\Sigma_g,\sigma)\ .$$
b) Suppose that the scalar field $\fk$ has characteristic zero,
 then the refined Verlinde formula is the following
\begin{eqnarray*} d_\mathcal{C}(\Sigma_g,\sigma)&=&
\langle \Omega\rangle^{g-1}
  \sum\limits_{\lambda\in\Gamma}\ 
  \la\lambda\ra^{2-2g} \times\prod\limits_{\nu=1}^g
\frac{\epsilon_\lambda(a_\nu(\sigma),b_\nu(\sigma))}
{(\sharp\, \mathrm{orb}(\lambda))^2}
\ .\end{eqnarray*}
\end{thm}
Here $(a(\sigma),b(\sigma))\in (\mathbb{Z}/d)^g\times (\mathbb{Z}/d)^g$
is given by the values of $q_\sigma$ on a sympleptic basis.

\esp\noindent If $\sharp\mathrm{orb}(\lambda)$ is even, then
$$\epsilon_\lambda(\mathfrak{a},\mathfrak{b})=\left\{\begin{array}{l}
1\ \text{ if  $\mathfrak{a}$ and $\mathfrak{b}$ are zero modulo $|Stab(\lambda)|$,}\\
0 \ \text{ else.}\end{array}\right.$$
If $\sharp\mathrm{orb}(\lambda)$ is odd, then
$$\epsilon_\lambda(\mathfrak{a},\mathfrak{b})=\left\{\begin{array}{l}
\frac{1}{2}(-1)^{\frac{2\mathfrak{a}}{|Stab(\lambda)|}\frac{2\mathfrak{b}}{|Stab(\lambda)|}}
\ \text{ if  $\mathfrak{a}$ and $\mathfrak{b}$ 
are  zero modulo $\frac{|Stab(\lambda)|}{2}$,}\\
0\ \text{ else.}\\
\end{array}\right.$$
\begin{rem}
Any element in $Stab(\lambda)$ has quantum dimension equal to one. In the case
where $\mathfrak{d}=-1$, the group $Stab(\lambda)$ is generated by an even power
 of $\varrho$, and $\sharp \mathrm{orb}(\lambda)$
is even.  We do not know examples with $\mathfrak{d}=-1$.
\end{rem}\begin{rem}
 If the scalar field $\fk$ has  characteristic $p>0$, then
statement b) computes the dimension mod. $p$.
\end{rem}
\begin{proof}
The formula in a) follows from the decomposition of the vector space
$\V_{\mathcal{C}}(\Sigma_g)$ described above.
 Moreover the dimension $d_\mathcal{C}(\Sigma_g,\sigma)$
 of a summand is the trace
of the corresponding projector. This projector can be represented
by a cobordism $ [0,1]\times\Sigma_g$ in which we have inserted a convenient
skein element. 
 By a standard TQFT argument we get
$$d_\mathcal{C}(\Sigma_g,\sigma)
=Z_\mathcal{C}(S^1\times \Sigma_g,\text{skein element})\ .$$
The $3$-manifold $S^1\times \Sigma_g$ is obtained by surgery on the borromean
link with $2g+1$ components represented in Figure \ref{cborr} \cite[Th. 14.12]{Li}.

 In this presentation, a meridian around the bigger 
component corresponds to $S^1\times pt$, and the $2g$ meridians around the
other components correspond to a system of $2g$ fundamental curves in $\Sigma_g$; these curves are framed by using the meridian disc.
 The skein element which arises here is represented by these
$2g$ curves, decorated with some $\pi_\nu$.
 
If $\mathfrak{d}=1$, then $\nu$ is the value of
the quadratic form $q_\sigma$ on the curve, and
if $\mathfrak{d}=-1$, then $\nu$ is the value of $\sigma$
on the $1$-cycle represented by the curve.

By using (\ref{normalized}) and lemma (\ref{proj}) we get
\begin{equation}d_\mathcal{C}(\Sigma_g,\sigma)=\mathcal{D}^{-(2g+2)}\sum_\lambda \langle \lambda \rangle
B_\lambda
=\langle \Omega\rangle^{-1-g}
\sum_\lambda \langle \lambda \rangle
B_\lambda\ ,\end{equation}
where $B_\lambda$ is the invariant of the colored borromean link in figure
\ref{cborr}. Here $(a_1,b_1),\dots,(a_g,b_g)$ are given by the values
of the quadratic form $q_\sigma$ on the corresponding curves
if $\mathfrak{d}=1$, and are equal to the value of $\sigma$
on the $1$-cycle represented by the curve
if $\mathfrak{d}=-1$.

\begin{figure}
\centerline{
\begin{picture}(40,65)
\put(0,5){\includegraphics[height=60mm]{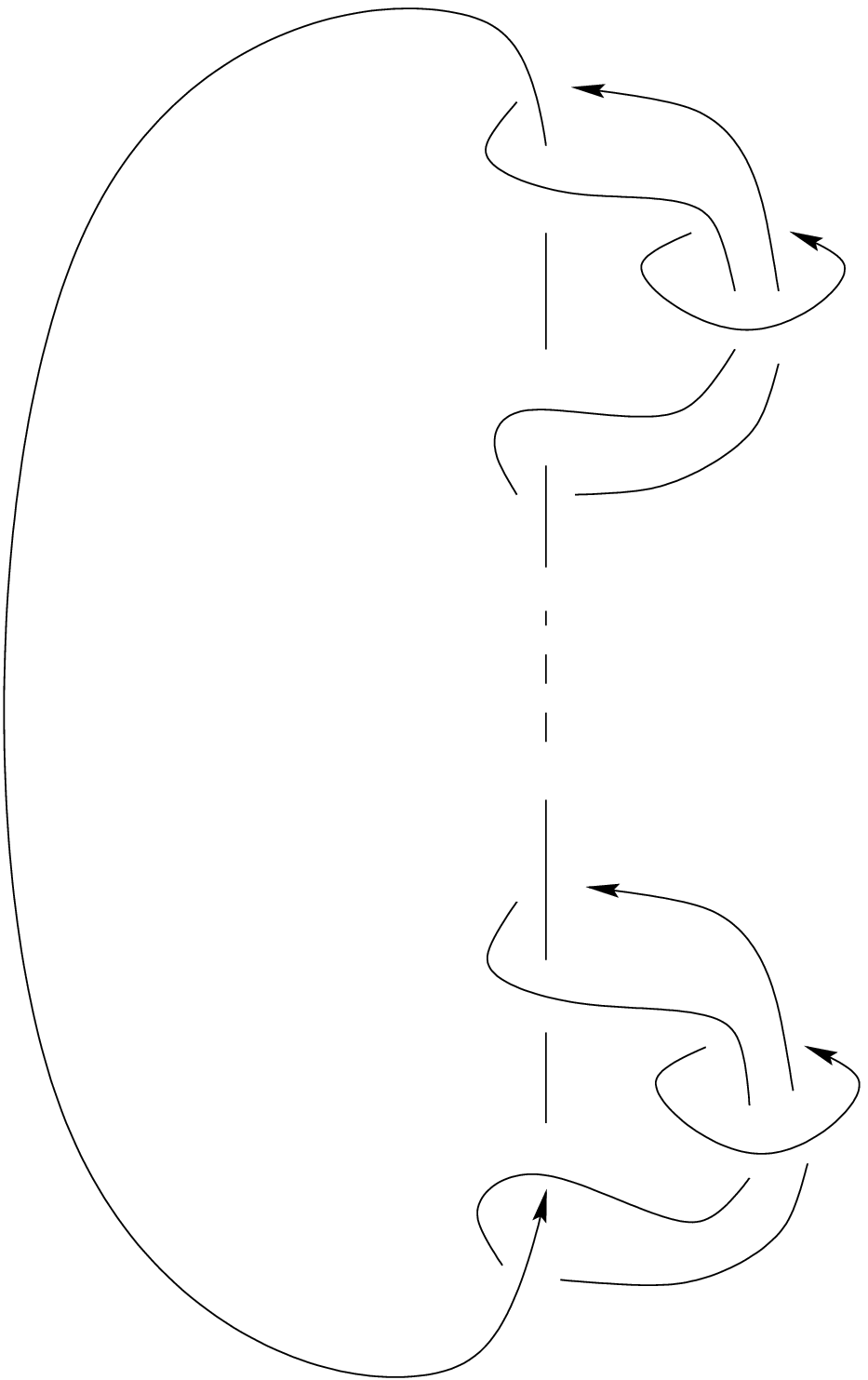}}
\put(1,10){$\lambda$}
\put(37,54){$\Omega_{a_1}$}\put(38,18){$\Omega_{a_g}$}
\put(32,60){$\Omega_{b_1}$}\put(33,25){$\Omega_{b_g}$}
\end{picture}
}
\caption{\label{cborr} Colored borromean link}
\end{figure}
Recall that the cyclic group generated by the
object $\varrho$, identified with $\mathbb{Z}/d$, acts on the set $\Gamma$ of
(representatives of) isomorphism classes of simple objects.
 If $j$ is in the stabilizer subgroup of $\lambda$,
then we choose a basis for the $1$-dimensional vector spaces
 $Hom_\mathcal{C}(\varrho^j,\lambda^*\otimes\lambda)$
and the dual basis for $Hom_\mathcal{C}(\lambda^*\otimes\lambda,\varrho^j)$.
 We denote these bases by the trivalent vertices in Figure \ref{tri}.

\begin{figure}\centerline{
\begin{picture}(30,30)
\put(3,3){\includegraphics[height=24mm]{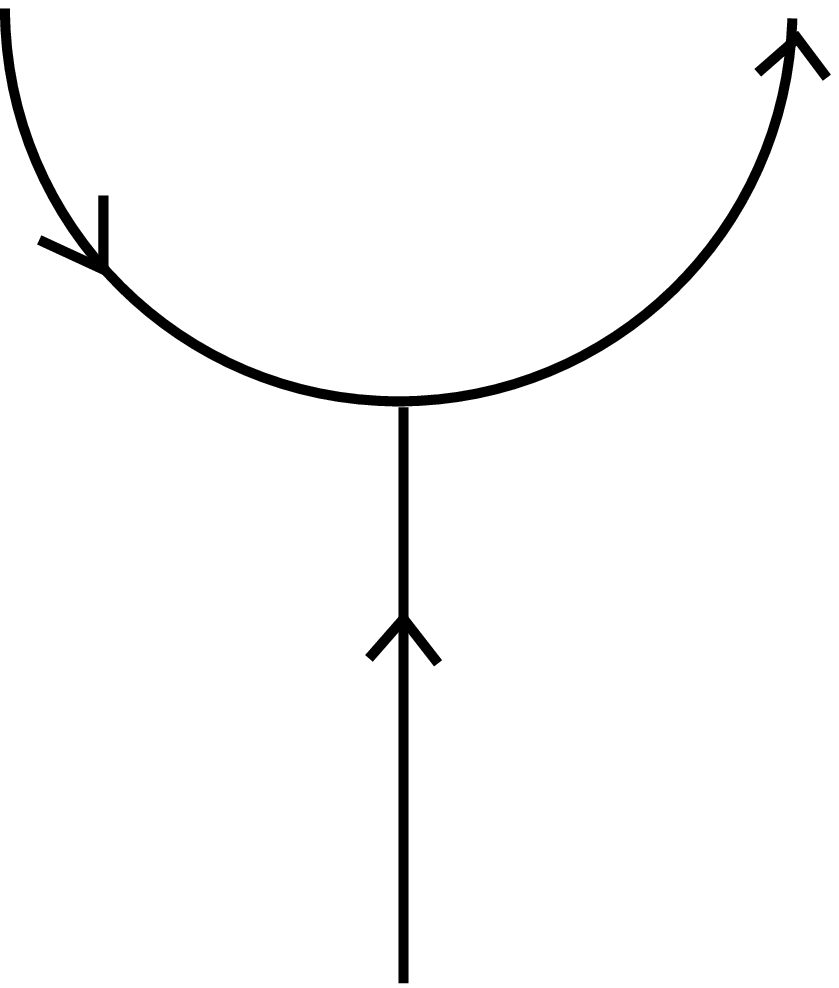}}
\put(1,20){$\lambda$}
\put(22,20){$\lambda$}
\put(14,8){$\varrho^j$}
\end{picture}
\hspace{2cm} 
\begin{picture}(30,30)
\put(3,3){\includegraphics[height=24mm]{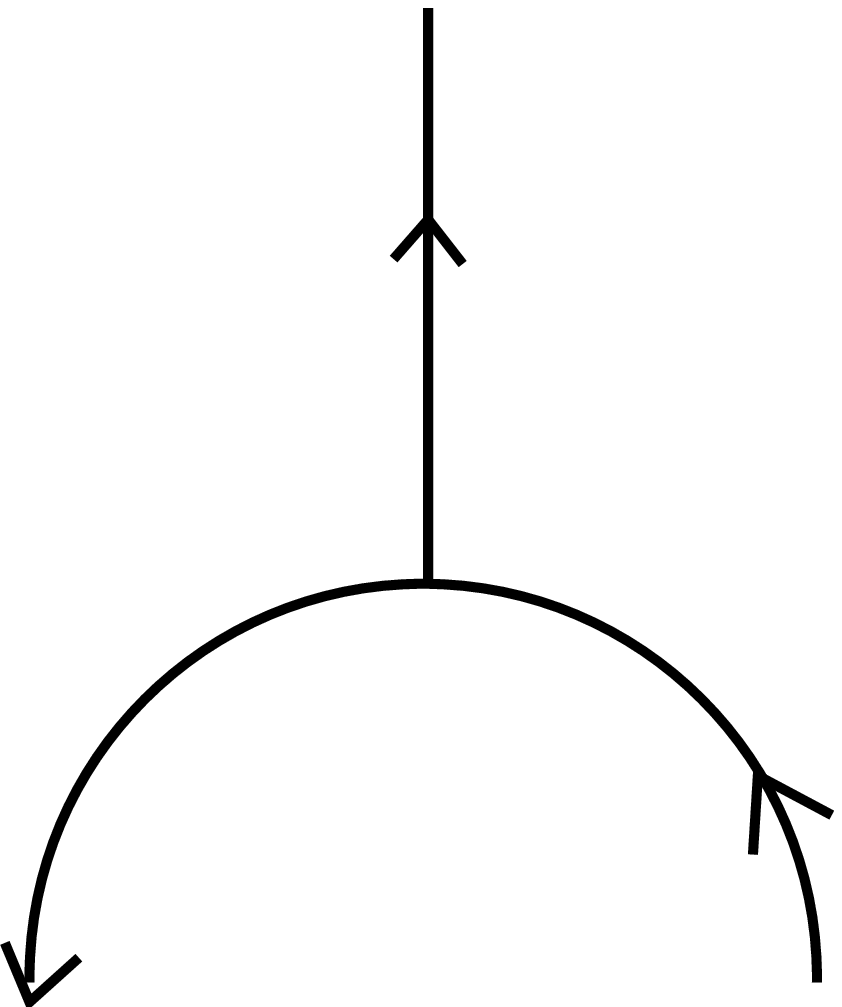}}
\put(1,7){$\lambda$}
\put(24,7){$\lambda$}
\put(16,18){$\varrho^j$}
\end{picture}
}
\caption{\label{tri}Trivalent vertices.}
\end{figure}

\begin{figure}
\centerline{
\begin{picture}(30,50)
\put(3,3){\includegraphics[height=40mm]{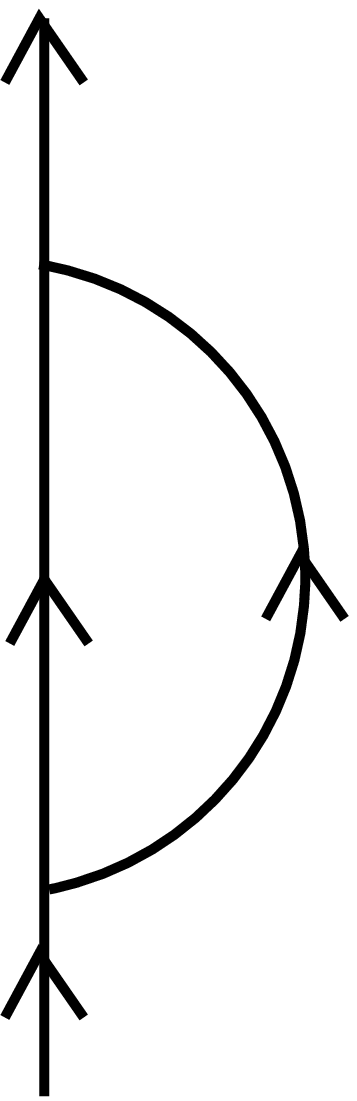}}
\put(0,25){$\lambda$}
\put(7,6){$\lambda$}\put(7,37){$\lambda$}
\put(15,25){$\varrho^j$}
\end{picture}
\raisebox{20mm}{$\displaystyle \ =\frac{1}{<\lambda>}\ $}
\begin{picture}(30,50)
\put(5,3){\includegraphics[height=40mm]{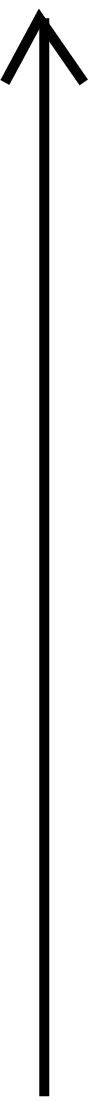}}
\put(8,27){$\lambda$}
\end{picture}}

\vspace{10pt}

\centerline{
\begin{picture}(30,50)
\put(3,3){\includegraphics[height=40mm]{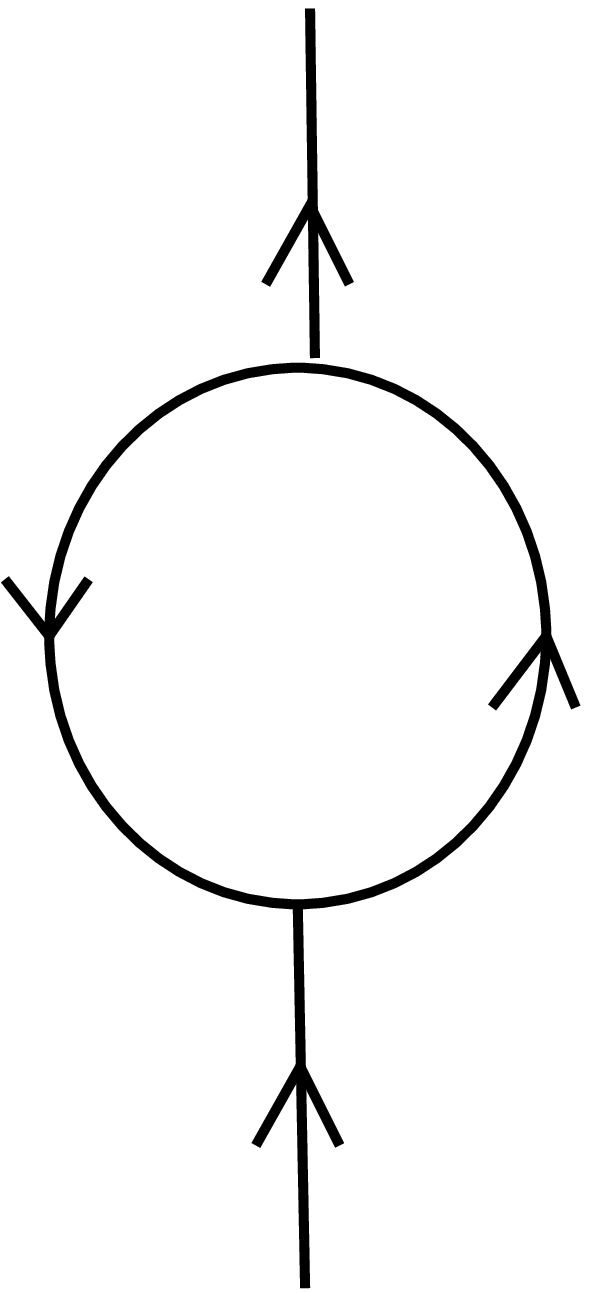}}
\put(1,27){$\lambda$}
\put(7,6){$\varrho^j$}\put(7,37){$\varrho^j$}
\put(20,27){$\lambda$}
\end{picture}
\raisebox{20mm}{$\displaystyle \ =\ \ $}\begin{picture}(30,50)
\put(5,3){\includegraphics[height=40mm]{rel3.eps}}
\put(8,27){$\varrho^j$}
\end{picture}}
\caption{\label{rel23}Relations for trivalent vertices.}

\end{figure}
We then have the  relations in Figure \ref{rel23};
recall that in the case $\mathfrak{d}=-1$, $j$ must be even. 

The proposition \ref{keypro} below is the key point in the computation.
 By using this proposition, we get
\begin{equation}\label{blambda}B_\lambda=\langle \lambda\rangle
\prod\limits_{\nu=1}^g
\sum_{j,j'\in Stab(\lambda)} \zeta^{ja_\nu+j'b_\nu}\frac{(-1)^{jj'}}
{\langle \lambda\rangle^2}
\frac{\langle \Omega\rangle^2}{d^2}
\end{equation}
Let $l=\sharp\mathrm{orb}(\lambda)$ and $l'=\frac{d}{l}=|Stab(\lambda)|$. The stabilizer subgroup
is then $Stab(\lambda)=\{ ls, \ 0\leq s <l'\}$.
 
If $l$ is even then $\sum_{j,j'\in Stab(\lambda)} \zeta^{ja_\nu+j'b_\nu}$
 is zero unless $\zeta^{la_\nu}=\zeta^{lb_\nu}=1$, and we get
$$\sum_{j,j'\in Stab(\lambda)} \zeta^{j a_\nu+j'b_\nu}(-1)^{jj'}=
\left\{
\begin{array}{ll}
|Stab(\lambda)|^2\ &\text{if $a_\nu\equiv b_\nu\equiv 0$ mod. $|Stab(\lambda)|$,}\\
0& \text{else.}
\end{array}
\right.$$

If $l$ is odd, then we decompose the sum 
$\sum_{j,j'\in Stab(\lambda)} \zeta^{ja_\nu+j'b_\nu}(-1)^{jj'}$ according to the parity of the indices.
 The sum  is zero if $\zeta^{2la_\nu}\neq 1$ or $\zeta^{2lb_\nu}\neq1$.
It remains four cases to consider according to  $\zeta^{la_\nu}=\pm 1$, $\zeta^{lb_\nu}=\pm 1$. 

Case $a_\nu\equiv b_\nu\equiv 0$ mod. $|Stab(\lambda)|$.
$$\begin{array}{ll}
\displaystyle\sum_{j,j'\in Stab(\lambda)} \zeta^{j a_\nu+j'b_\nu}(-1)^{jj'}&=\displaystyle\sum_{j,j'\text{ even}}
+\sum_{j\text{ even ,}j' \text{ odd}}+\sum_{j'\text{ even ,}j \text{ odd}}-\sum_{j,j' \text{ odd}}\\
&=\displaystyle\frac{{l'}^2}{4}+\frac{{l'}^2}{4}+\frac{{l'}^2}{4}-\frac{{l'}^2}{4}\\
&=\displaystyle\frac{{l'}^2}{2}\ .\end{array}$$
Case $a_\nu\equiv b_\nu\equiv  \frac{|Stab(\lambda)|}{2}$ mod. $|Stab(\lambda)|$.
$$\begin{array}{ll}
\displaystyle\sum_{j,j'\in Stab(\lambda)} \zeta^{j a_\nu+j'b_\nu}(-1)^{jj'}&=\displaystyle\sum_{j,j'\text{ even}}
+\sum_{j\text{ even ,}j' \text{ odd}}+\sum_{j'\text{ even ,}j \text{ odd}}-\sum_{j,j' \text{ odd}}\\
&=\displaystyle\frac{{l'}^2}{4}-\frac{{l'}^2}{4}-\frac{{l'}^2}{4}-\frac{{l'}^2}{4}\\
&=-\displaystyle\frac{{l'}^2}{2}\ .\end{array}$$
Case $a_\nu\equiv 0$ , $b_\nu\equiv  \frac{|Stab(\lambda)|}{2}$ mod. $|Stab(\lambda)|$.
$$\begin{array}{ll}
\displaystyle\sum_{j,j'\in Stab(\lambda)} \zeta^{j a_\nu+j'b_\nu}(-1)^{jj'}&=\displaystyle\sum_{j,j'\text{ even}}
+\sum_{j\text{ even ,}j' \text{ odd}}+\sum_{j'\text{ even ,}j \text{ odd}}-\sum_{j,j' \text{ odd}}\\
&=\displaystyle\frac{{l'}^2}{4}+\frac{{l'}^2}{4}-\frac{{l'}^2}{4}+\frac{{l'}^2}{4}\\
&=\displaystyle\frac{{l'}^2}{2}\ .\end{array}$$
Case $b_\nu\equiv 0$ , $a_\nu\equiv  \frac{|Stab(\lambda)|}{2}$ mod. $|Stab(\lambda)|$.
$$\begin{array}{ll}
\displaystyle\sum_{j,j'\in Stab(\lambda)} \zeta^{j a_\nu+j'b_\nu}(-1)^{jj'}&=\displaystyle\sum_{j,j'\text{ even}}
+\sum_{j\text{ even ,}j' \text{ odd}}+\sum_{j'\text{ even ,}j \text{ odd}}-\sum_{j,j' \text{ odd}}\\
&=\displaystyle\frac{{l'}^2}{4}-\frac{{l'}^2}{4}+\frac{{l'}^2}{4}+\frac{{l'}^2}{4}\\
&=\displaystyle\frac{{l'}^2}{2}\ .\end{array}$$

In all cases we get the formula below.
\begin{equation}\label{sum}
\sum_{j,j'\in Stab(\lambda)} \zeta^{j a_\nu+j'b_\nu}(-1)^{jj'}=
\epsilon_\lambda(a_\nu,b_\nu)\ |Stab(\lambda)|^2\ .
\end{equation}
 
In the case where $\mathfrak{d}=-1$, then  $l$ is even, $l'=|Stab(\lambda)|$ divides $\frac{d}{2}$ and we have
$\epsilon_\lambda(a_\nu,b_\nu)=\epsilon_\lambda(a_\nu+\frac{d}{2},b_\nu+\frac{d}{2})$.
 So that we may define $(a_1,b_1),\dots,(a_g,b_g)$  by the values
of the quadratic form $q_\sigma$ as well.

We will now establish statement b). 
$$\begin{array}{ll}d_\mathcal{C}(\Sigma_g,\sigma)&\displaystyle
=\langle \Omega\rangle^{-1-g}
\sum_\lambda \langle \lambda \rangle^{2-2g}
\prod_{\nu=1}^g \epsilon_\lambda(a_\nu,b_\nu)\ |Stab(\lambda)|^2\ \frac{\langle \Omega\rangle^2}{d^2}
\\&\displaystyle
=\langle \Omega\rangle^{g-1}
\sum_\lambda \langle \lambda \rangle^{2-2g}
\prod_{\nu=1}^g \frac{\epsilon_\lambda(a_\nu,b_\nu)}{\sharp \mathrm{orb}(\lambda)^2}\ .
\end{array}$$
\end{proof}
\begin{lem}\label{exists}
 For any $i\in\mathbb{Z}/d$ the subcategory $\mathcal{C}_i$ contains at least
one simple object, and for any simple object  $\lambda_i$ in $\mathcal{C}_i$, 
one has 
$$\langle \lambda_i\rangle \Omega_{i+j}=\lambda_i\otimes \Omega_j\ .$$
\end{lem}
In a modular category the dimension of a simple object is non-zero,
hence we have that  for any $i$
\begin{equation}\langle \Omega_i\rangle =\langle \Omega_0\rangle=
\frac{1}{d}\langle \Omega\rangle\end{equation}
\begin{proof}
Let $\nu$ be a generator for the subgroup
of $\mathbb{Z}/d$ formed with all $i$ such that $\mathcal{C}_i$
contains at least one non-trivial object.
 Suppose that $\nu$ has order $d'$, then $\varrho^{d'}$
is a simple object whose contribution in the $S$ matrix is the
same as that of the trivial. This object is isomorphic to the trivial,
and we deduce that $d'=d$.
 This proves the first part of the lemma. The second part follows from
the graded sliding property (see \cite[Section 1.3]{LT}).
\end{proof}
\begin{lem}\label{vanish}
 Let $\lambda$ be a simple object in $\mathcal{C}$, for any $i$ in $\mathbb{Z}/d$   
 the following morphism 
is non-zero   if and only if $\lambda$ is isomorphic to
$\varrho^j$ for some $j$.

\centerline{\includegraphics[height=25mm]{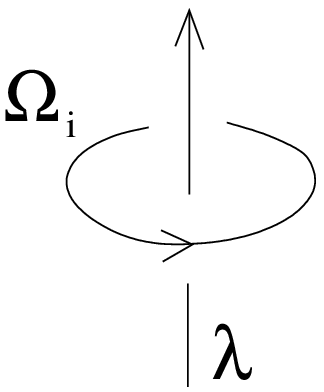}}
\end{lem}
\begin{proof}
If $\lambda$ is equal to $\varrho^j$ then the morphism
is equal to $\frac{1}{d}\langle\Omega\rangle \zeta^{ij}\one\!_\lambda$,
and so is not zero.
 
 Suppose now that for some simple object $\lambda$
the above morphism is not zero.
 By using Lemma \ref{exists}, we obtain a scalar $t_\lambda$ such that

\vspace{5pt}
\centerline{\includegraphics[height=25mm]{circbcd.eps}
\raisebox{11mm}{$=
t_\lambda^i$}\includegraphics[height=25mm]{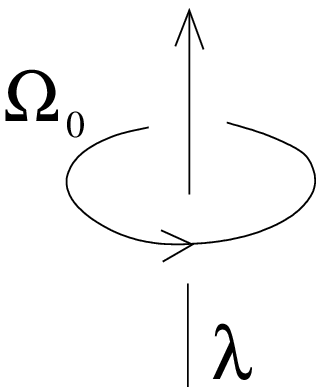}\ .}
\vspace{5pt}

Note that $t_\lambda^d=1$, hence there exists $j$ such that
 $t_\lambda=\zeta^j$. By the graded sliding property
we deduce that the contribution of $\lambda$ in the $S$ matrix is the same
as that of $\varrho^j$, and we get the required isomorphism.
\end{proof}
\begin{lem}
\label{rfusion}
 Let $\lambda$ be a simple object in $\mathcal{C}$, then one has the relation
in Figure \ref{rfus}.
\begin{figure}
\centerline{
\begin{picture}(30,50)
\put(3,5){\includegraphics[height=40mm]{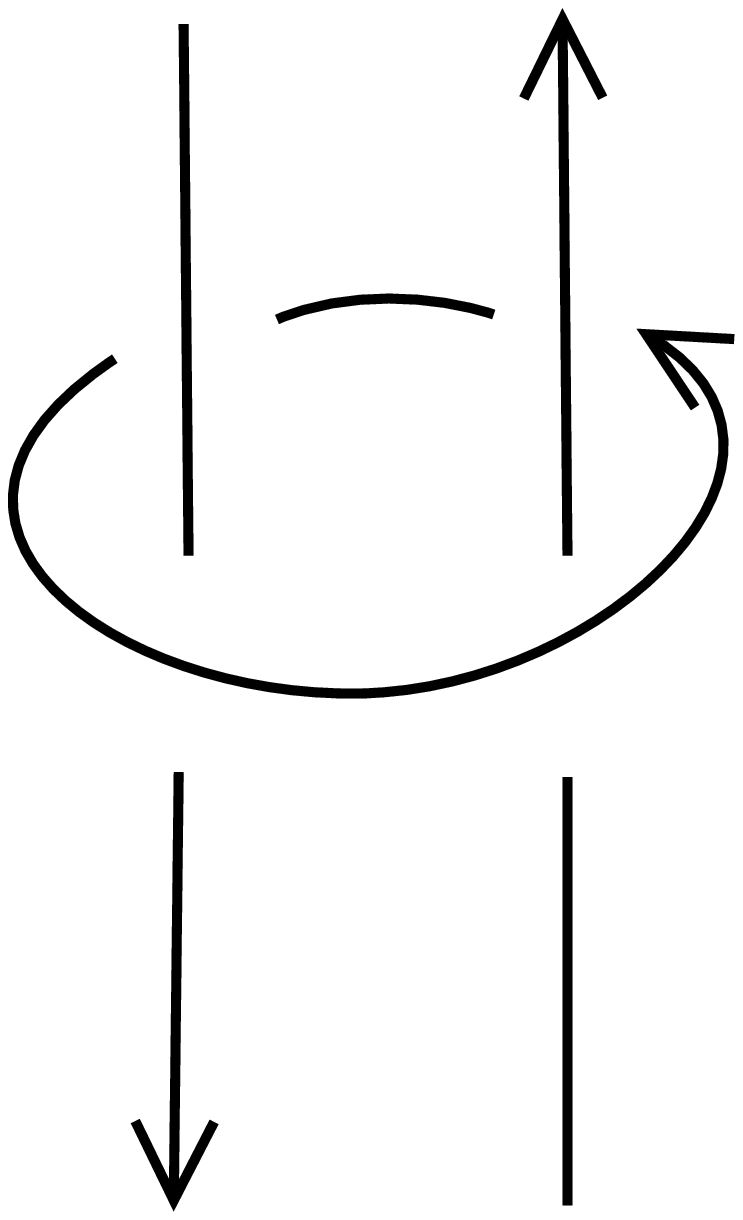}}
\put(4,27){$\Omega_a$}
\put(11,40){$\lambda$}
\put(28,40){$\lambda$}
\end{picture}
\raisebox{22mm}{$\ \ \ =\ \ \sum_{j\in Stab(\lambda)}
\zeta^{aj}\ \frac{\langle \Omega \rangle}{d}\ $}
\begin{picture}(30,50)
\put(4,5){\includegraphics[height=40mm]{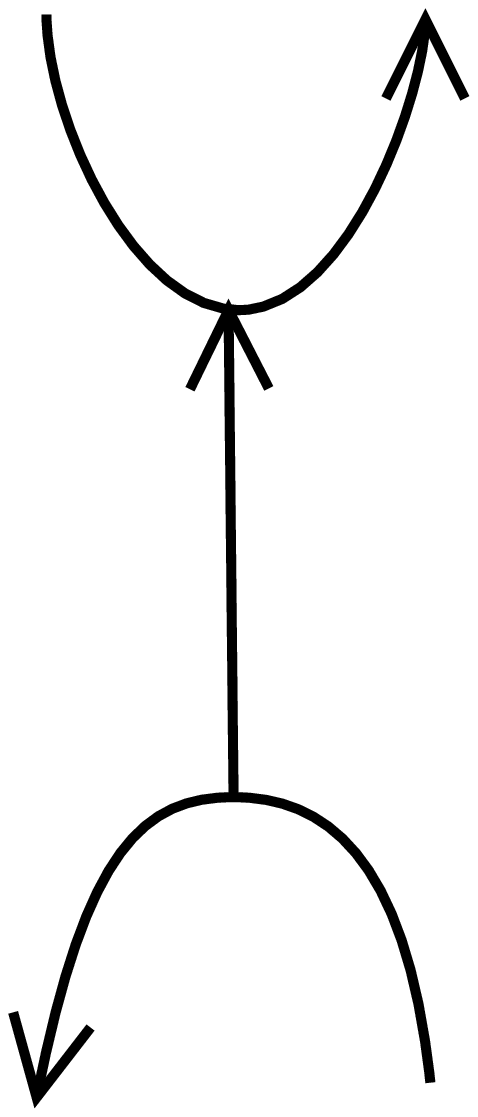}}
\put(14,27){$\varrho^j$}
\put(3,40){$\lambda$}
\put(21,40){$\lambda$}
\put(3,7){$\lambda$}
\put(20,7){$\lambda$}
\end{picture}
}
\caption{\label{rfus}}
\end{figure}
\end{lem}
\begin{proof}
We first use the domination axiom. The decomposition of the identity
of $\lambda^*\otimes\lambda$ is given by a so called {\em fusion formula}
(see e.g. \cite[Section 1.2]{BB}). Note that in this formula the multiplicity
of an invertible objet is one if it belongs to the stabilizer subgroup
of $\lambda$ and zero else.
 We then apply Lemma \ref{vanish}.
 The result follows.
\end{proof}
\begin{lem}\label{tetra} For $i,j$ in $Stab(\lambda)$, one has
the relation in Figure \ref{tet}.
\end{lem}
\begin{figure}
\centerline{
\begin{picture}(25,50)
\put(10,5){\includegraphics[height=40mm]{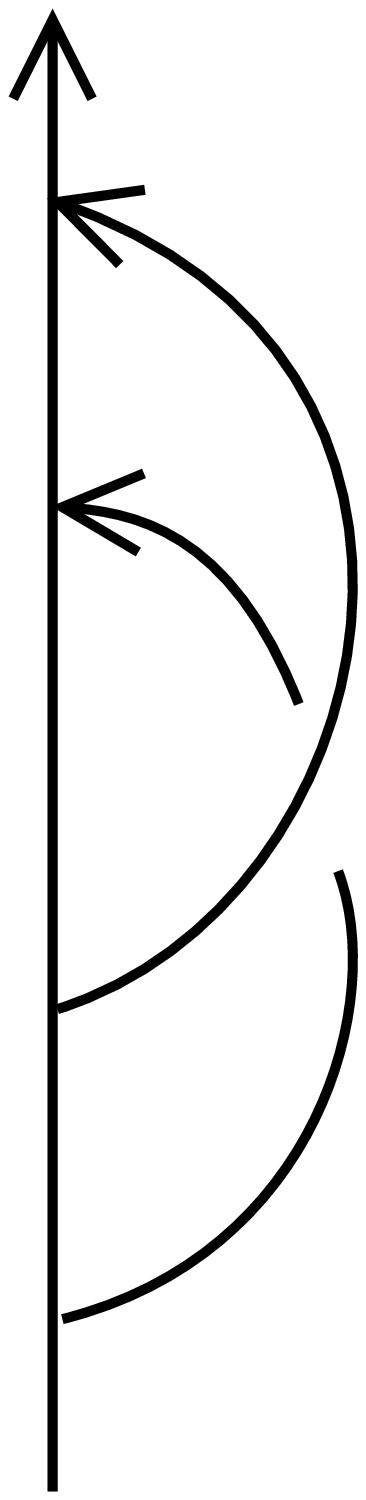}}
\put(8.5,15){$\lambda$}
\put(19,34){$\varrho^i$}
\put(18,12){$\varrho^j$}
\end{picture}
\raisebox{22mm}{$
 \ =\ {(-1)^{ij}}\ $
} 
\begin{picture}(18,50)
\put(3,5){\includegraphics[height=40mm]{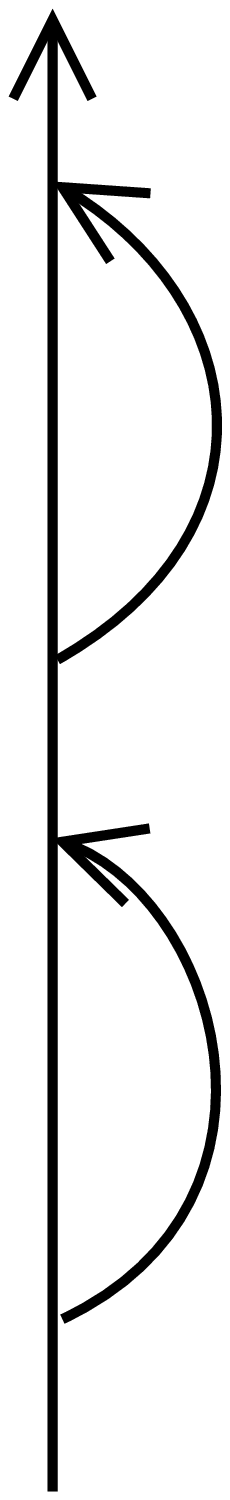}}
\put(1,15){$\lambda$}
\put(11,34){$\varrho^i$}
\put(10,13){$\varrho^j$}
\end{picture}
\raisebox{22mm}{$
\  =\ \frac{(-1)^{ij}}{\langle \lambda\rangle^2}\ $
} 
\begin{picture}(20,50)
\put(3,5){\includegraphics[height=35mm]{rel3.eps}}
\put(7,25){$\lambda$}
\end{picture}
}
\caption{\label{tet}}
\end{figure}
\begin{proof}
The first equality uses the defining property of a 
modulo $d$ spin modular category.
 The second one comes from the definition of the trivalent vertices.
\end{proof}
\begin{pro}\label{keypro}
The  formula in Figure \ref{key} holds.
\begin{figure}
\centerline{\begin{picture}(30,50)
\put(3,5){\includegraphics[height=40mm]{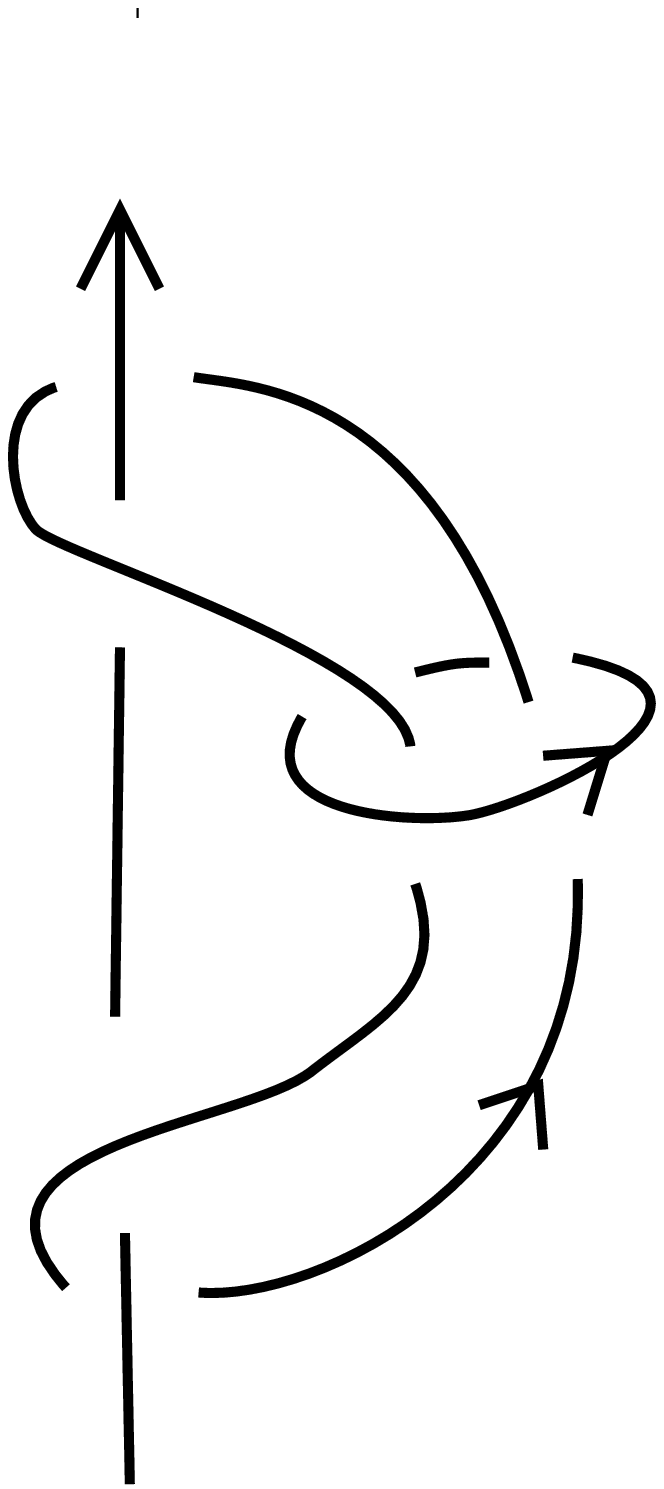}}
\put(3,20){$\lambda$}
\put(16,34){$\Omega_b$}
\put(23,27){$\Omega_a$}
\end{picture}
\raisebox{22mm}{$\ =\ \sum_{j,j'\in Stab(\lambda)}\ \
\zeta^{aj+bj'}\ \frac{(-1)^{jj'}}
{\langle \lambda\rangle^2}
\ \frac{\langle \Omega \rangle^2}{d^2}
\ $}
\begin{picture}(30,50)
\put(3,5){\includegraphics[height=35mm]{rel3.eps}}
\put(7,20){$\lambda$}
\end{picture}
}
\caption{\label{key}}
\end{figure}
\end{pro}
\begin{proof}
By using Lemma \ref{rfusion} twice (firstly for the component colored by $\Omega_b$), we get
the formula in Figure \ref{keycalc}.
\begin{figure}
\centerline{
\raisebox{22mm}{$lhs\  =\ \sum_{j,j'\in Stab(\lambda)}
\zeta^{aj+bj'}\ \frac{\langle \Omega \rangle^2}{d^2}\ $}
\begin{picture}(50,50)
\put(3,5){\includegraphics[height=40mm]{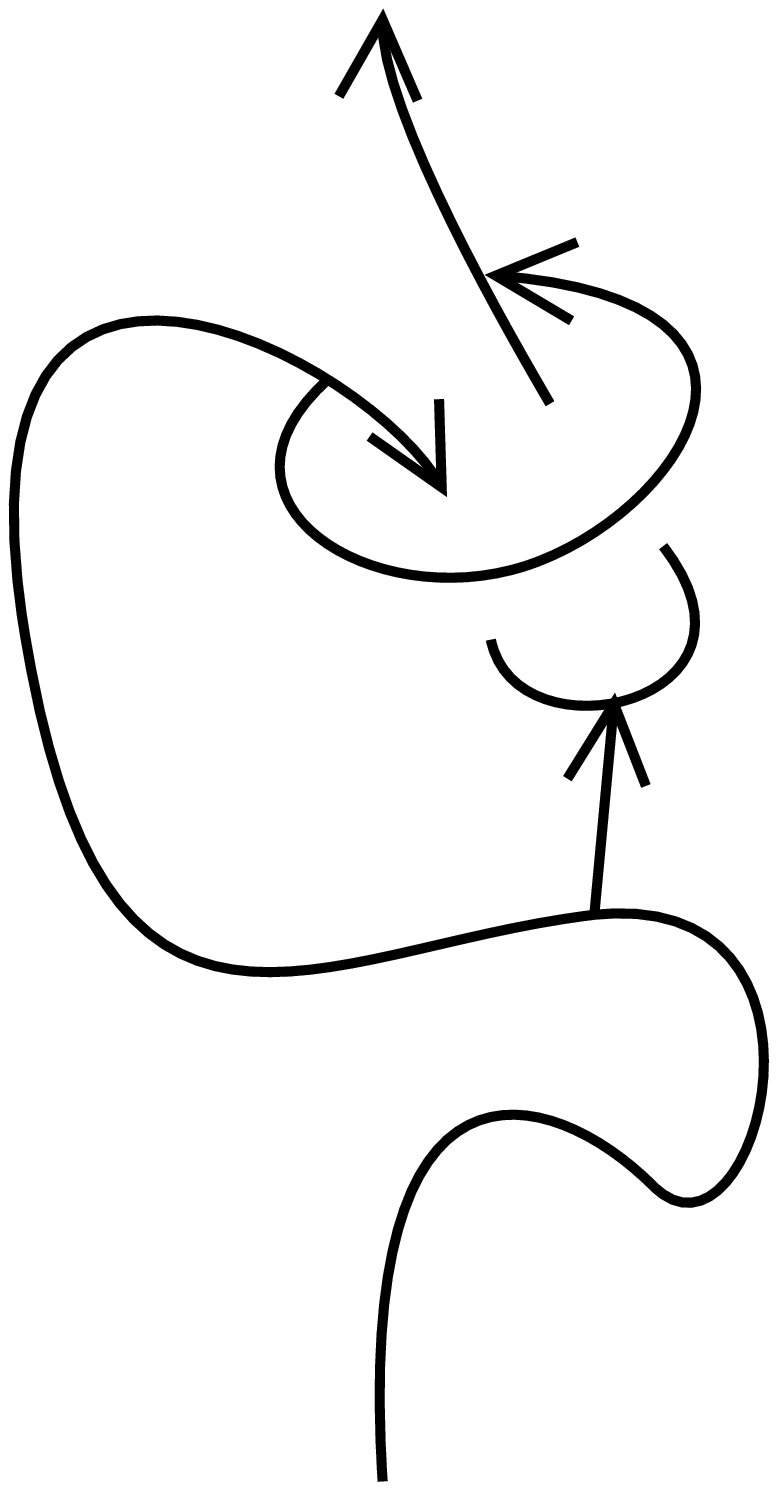}}
\put(24,35){$\varrho^{j'}$}
\put(23,22){$\varrho^{j}$}
\put(17,10){$\lambda$}
\end{picture}
}
\caption{\label{keycalc}}
\end{figure}
After an isotopy, we  apply Lemma \ref{tetra}. The result follows.
\end{proof}

\section{Verlinde formulas for type A modular categories}
\subsection{The $SU(N,K)$ modular category}
We first consider the so called $SU(N,K)$ modular category.
 The construction can be done either from the representation theory
of the quantum group $U_qsl(N)$ at a convenient root of unity
\cite{AP,TW,BK} or from skein theory \cite{Yo,Bhec}.
 In the following we will use Young diagrams to denote the corresponding simple object.
 Here a Young diagram (or partition)
$\lambda$ is a finite non-increasing sequence
of non-negative integers.
A cell for this partition is a pair $c=(i,j)$
with $1\leq j\leq \lambda_i$.
 We denote by $\lambda^{\vee}$ the tranpose of $\lambda$;
$(i,j)$ is a cell in  $\lambda^{\vee}$ if and only if $(j,i)$
is a cell in  $\lambda$.
 The content and hook-length for a cell $c=(i,j)$
are defined respectively by
$$cn(c)=j-i\ , \ hl(c)=\lambda_i+\lambda^{\vee}_j-i-j+1\ .$$
The size of $\lambda$ is $|\lambda|=\sum_i \lambda_i$.

The following theorem is proved in \cite{Bhec}. The result can also be obtained
from \cite[Th. 3.3.20]{BK} ($A_{N-1}$ case).
\begin{thm} \label{modcatH}
Let $N,K\geq 2$. Suppose that $a$ is a $2N(N+K)$-th root
of unity in the scalar field and $s=a^{-N}$.

There exists a modular category
 $\mathcal{C}^{SU(N,K)}$ whose set of distinguished  simple objects
 is 
 $$\Gamma_{N,K}=\{\lambda=(\lambda_1,\dots,\lambda_{N}),\
K\geq\lambda_1\geq\dots\geq \lambda_{N-1}\geq\lambda_N=0\}\ .$$

The quantum dimension and framing coefficient
of a simple object $\lambda\in \Gamma_{N,K}$ are given by the following formulas
(here $[n]=\frac{s^n-s^{-n}}{s-s^{-1}}$ denotes the quantum integer).

$$\langle \lambda \rangle=
\prod_{\mathrm{cells}}\frac{[N+cn(c)]}{[hl(c)]}$$

\centerline{\begin{pspicture}[.5](0,-0.5)(1,1)
\pscurve{-}(.1,0)(.15,.45)(.45,.75)(.7,.5)(.4,.2)(.20,.38)
\pscurve{->}(.14,.62)(.1,.9)(.1,1)
 \put(-.2,0){$\lambda$}
\end{pspicture}
$\ =\ a^{|\lambda|^2}s^{N|\lambda|+2\sum_{\mathrm{cells}}cn(c)}\ \ $
\begin{pspicture}[.5](0,-.5)(1,1)
\psline{->}(.1,0)(.1,1)
 \put(.2,0){$\lambda$}
\end{pspicture}
}
\end{thm}
\begin{rems}1. In the quantum group approach, a Young diagram in $\Gamma_{N,K}$
gives a highest weight module, which is irreducible and has non zero quantum dimension. 
 The quantum dimension follows from Weyl's character formula and computation with symmetric functions in \cite[section I.3]{Macdo}. The value of the twist is obtained by the action of Drinfeld quantum Casimir.\\
2. In the skein theoretic approach the Young diagram gives a minimal idempotent in Hecke algebra (the deformation of the Young symmetrizer
in the symmetric group algebra); this idempotent becomes a simple object
in the so called Karoubi completion of the Hecke category.
\end{rems}
We denote by $\V_{N,K}(\Sigma_g)$
the TQFT vector space, associated with a genus $g$ surface
$\Sigma_g$,
 for the modular category $\mathcal{C}^{SU(N,K)}$, and by
$d_{N,K}(g)$
 its rank. We give below the well known computation for this formula.

\begin{thm}
The rank $d_{N,K}(g)$ is
equal to the Verlinde number for the group
$SU(N)$ at level $K$.
$$\begin{array}{l}
d_{N,K}(g)=\mathcal{V}_{SU(N)}(K,g)=\hfill \\
  \left((N+K)^{(N-1)}N\right)^{g-1}
  \sum\limits_{\lambda\in\Gamma_{N,K}}\ \prod\limits_{1\leq i<j\leq N}
  \left({2\sin{(\lambda_i-i-\lambda_j+j)\frac{\pi}{N+K}}}
\right)^{2-2g} 
  \end{array}$$
\end{thm}
\begin{proof}
We can use Turaev's formula \cite[Corollary 12.1.2]{Tu}.
Note that this formula computes the TQFT-invariant of the manifold
$S^1\times \Sigma_g$.
$$d_{N,K}(g)=\left(\sum_{\lambda\in\Gamma_{N,K}}\langle \lambda\rangle^2\right)
^{g-1}\sum\limits_{\lambda\in\Gamma_{N,K}}
\langle \lambda \rangle^{2-2g}$$
The computation  is achieved with the lemma below.
Statement a) is a standard fact on symmetric functions; statement b)
is contained e.g. in  the proof of lemma 2.8
in \cite{Bhec}.
 Note that the result does not depend on the choice
of the root of unity with required order; it is also unchanged if $s$ is replaced by $\tilde s=-s$.
 The formula agrees with the Verlinde
number $\mathcal{V}_{SU(N)}(K,g)$ \cite{Beau,So}.
\end{proof}
 \begin{lem}
a) $$\langle \lambda\rangle^2=
\frac{a_{\rho+\lambda}\overline a_{\rho+\lambda}}
{a_{\rho}\overline a_{\rho}}$$
with $\rho=(N-1,N-2,\dots,0)$ and, for
$l=(l_1,\dots,l_N)$
$$a_l=det\left(s^{2(i-1)l_j}\right)
_{1\leq i,j\leq N}$$
b) $$ \sum_{\lambda\in\Gamma_{N,K}}\langle \lambda\rangle^2
=\frac{N(N+K)^{N-1}}{a_{\rho}\overline a_{\rho}}.$$
\end{lem}
\subsection{Spin decomposition of the Verlinde formula for ${SU(N,K)}$ modular category.}
In the category $\mathcal{C}^{SU(N,K)}$, the object $(K)$ (a $K$ cells Young diagram with only 
one row) is an invertible object whose
order is $N$. It is a generator of the group of invertible
objects, and has quantum dimension $1$. Its framing coefficient is
equal to 
$$\theta_K=a^{K^2}s^{NK+K(K-1)}=(-a^{N+K})^K\ .$$ 
If $N=jl$, and $(-a^{N+K})^{Kj^2}=-1$, then the category 
$\mathcal{C}^{SU(N,K)}$ equipped with the invertible object $\varrho=(K)^{\otimes j}$ is a 
modulo $l$ spin modular category.

Recall that $a$ is a $2N(N+K)$-th root of unity. Let $d=gcd(N,K)$, $N=dN'$,
$K=dK'$. A convenient integer $j$ exists if and only if
\begin{verse}
either $d$ is even, $N'$ is odd and the exponent of $2$ in  $K'$ is even 
($K'=2^{2n}(2m+1)$),\\
or $d$ is odd and the exponent of $2$ in $N$ is an even positive number.
\end{verse}
We emphasize the simplest case in the theorem below.
\begin{thm}
If $N$ is even and $K'=\frac{K}{N}$ is an odd integer, then the category
$\mathcal{C}^{SU(N,K)}$ equipped with the invertible object $\varrho=(K)$
 is a modulo $N$ spin  modular category.
\end{thm}
 The following theorem is an application of \ref{spin}.
\begin{thm}\label{spin_A}
Suppose that  $N$ is even and $K'=\frac{K}{N}$ is an odd integer.
a) There exists a splitting of the Verlinde formula
$$d_{N,K}(g)=\sum\limits_{\sigma\in Spin(\Sigma_g,\mathbb{Z}/d)}
d_{N,K}(g,\sigma)\ .$$
b) The refined Verlinde formula is the following
\begin{eqnarray*} d_{N,K}(g,\sigma)&=&
\left((N+K)^{N-1}N\right)^{g-1}
  \sum\limits_{\lambda\in{\Gamma}_{N,K}}\ \prod\limits_{\nu=1}^g
\frac{\epsilon_\lambda(a_\nu(\sigma),b_\nu(\sigma))}
{(\sharp\, \mathrm{orb}(\lambda))^2}\\
& &\times\prod\limits_{1\leq i<j\leq N}
  \left({2\sin{(\lambda_i-i-\lambda_j+j)\frac{\pi}{N+K}}}
\right)^{2-2g} 
\ .\end{eqnarray*}
\end{thm}
Here $(a(\sigma),b(\sigma))\in \mathbb{Z}/N)^g\times (\mathbb{Z}/N)^g$
are given by the values of $q_\sigma$ on a sympleptic basis.
 We consider the action of $\mathbb{Z}/N$ on the set
$$\Gamma_{N,K}=\{\lambda=(\lambda_1,\dots,\lambda_{N}),\
K\geq\lambda_1\geq\dots\geq \lambda_{N-1}\geq\lambda_N=0\}\ ,$$ 
given for the generator 
 of the cyclic group
$\mathbb{Z}/N$ by
$$(\lambda_1,\dots,\lambda_{N-1},0)\longmapsto
(K,\lambda_1,\dots,\lambda_{N-1})-(\lambda_{N-1},\dots,\lambda_{N-1})\ .$$ 
 We denote by $\sharp\, \mathrm{orb}(\lambda)$
 the cardinality of the orbit of $\lambda$,
and by $Stab(\lambda)$ the stabilizer subgroup. 
 The numbers $\epsilon_\lambda(\mathfrak{a},\mathfrak{b}) \in
\{0,1,-\frac{1}{2},\frac{1}{2}\}$ are defined  as follows.

\esp\noindent If $\sharp\mathrm{orb}(\lambda)$ is even, then
$$\epsilon_\lambda(\mathfrak{a},\mathfrak{b})=\left\{\begin{array}{l}
1\ \text{ if  $\mathfrak{a}$ and $\mathfrak{b}$ are zero modulo $|Stab(\lambda)|$,}\\
0 \ \text{ else.}\end{array}\right.$$
If $\sharp\mathrm{orb}(\lambda)$ is odd, then
$$\epsilon_\lambda(\mathfrak{a},\mathfrak{b})=\left\{\begin{array}{l}
\frac{1}{2}(-1)^{\frac{2\mathfrak{a}}{|Stab(\lambda)|}\frac{2\mathfrak{b}}{|Stab(\lambda)|}}
\ \text{ if  $\mathfrak{a}$ and $\mathfrak{b}$ 
are  zero modulo $\frac{|Stab(\lambda)|}{2}$,}\\
0\ \text{ else.}\\
\end{array}\right.$$
\begin{rem}
In the general case, one can use the reduction formula \cite[Theorem 3.6]{Bhec}
in order to establish a tensor product decomposition of the $SU(N,K)$
TQFT functor $\V_{N,K}$
$$\V_{N,K}=\V^{U(1)}_{N'}\otimes \widetilde{\V}_{N,K}\  ,$$ 
where $\widetilde \V_{N,K}$ is the TQFT functor
associated with the modular category
$\mathcal{C}^{PU(N,K)}$ discussed below, and $\V^{U(1)}_{N'}$
(known as a $U(1)$ theory)
is associated with a modular category based on linking numbers. The latter
involves a root of unity $\eta$ whose order is $2N'$ (resp. $N'$) if $N'$
is even (resp. odd); when $N'$ is even with even exponent of $2$, then one
can find $j$ such that $j^2\equiv N'$ ($mod\ 2N'$), and the category
is modulo $j$ spin  modular.
\end{rem}

\subsection{The $PU(N,K)$ modular category}
The so called projective  $PSU(N,K)$ modular category was obtained for $N$ and $K$ coprime
by restricting to simple objects in the root lattice \cite{MW,Le,LT}.
 The modular category
$\mathcal{C}^{PU(N,K)}$ (denoted by $\widetilde{\mathrm{H}}^{N,K}$
in \cite{Bhec}) is a generalization
 to the case where
$N$ and $K$ are not required to be  coprime.

Let $N,K\geq 2$. We suppose that 
 in the scalar field
$$\begin{cases}
s \text{ has order $2(N+K)$ if $N+K$ is even,}\\
s \text{ has order $N+K$ if $N+K$ is odd.}
\end{cases}$$ 
Then the Hecke category completed with idempotents and quotiented with
negligible, which we denote by $H^{N,K}$ is semisimple. In addition to simple objects $\lambda\in \Gamma_{N,K}$ there is an invertible simple object $1^N$ and its tensor powers.
 The group of invertible objects is generated by $1^N$ and $(K)$
with the relation $(1^N)^{\otimes K}\approx (K)^{\otimes N}$.
In order to apply the modularization procedure, we have to know
which are the transparent simple objects \cite{Bru}. This depends
on the order $\alpha$ of $(a^Ns)^2$ and  the order $\beta$ of
$(a^Ks^{-1})^2$. The set of isomorphism classes of transparent
simple objects is then the group generated by 
$(1^N)^{\otimes \alpha}$ and $(K)^{\otimes \beta}$.
We choose the framing parameter $a$ in such a way that this group
of transparent objects is as big as possible, and that the modularization
criterion is satisfied.
\begin{thm} \label{reduced}
Set $d=\mathrm{gcd}(N,K)$, $N=dN'$, $K=dK'$, $d=\alpha \beta$ with
$\mathrm{gcd}(\alpha,K')=\mathrm{gcd}(\beta,N')=
\mathrm{gcd}(\alpha,\beta)=1$.\\
Suppose that $a$ satisfies the relations
$$(a^Ns)^\alpha=(-1)^{N+K+1}(a^Ks^{-1})^\beta=(-1)^{(N+K+1)\beta}$$
(such an $a$ exists up to extension of the scalar field).

There exists a modular category
 $\mathcal{C}^{PU(N,K)}$ in which  isomorphism classes of
simple objects
corresponds bijectively with cosets in the quotient of 
$$\dot \Gamma_{N,K}=\{(1^N)^{\otimes j}\otimes \lambda,\ 0\leq j<\alpha,\
\lambda\in \Gamma_{N,K}\}$$
under a free action of the cyclic group of order $N/\alpha$.
\end{thm}
The action of the generator is given by tensor product with $(K)^{\otimes \beta}$
in the completed Homfly category $H$. One has to iterate $\beta$ times
the rule
$$(1^N)^{\otimes j}\otimes\lambda\mapsto (1^N)^{\otimes j'}+(K-\lambda_{N-1},
\lambda_1-\lambda_{N-1},\dots,\lambda_{N-2}-\lambda_{N-1},0)$$
where $j'\equiv j+\lambda_{N-1}\ \mathrm{mod}\ \alpha$.

The quantum dimension and framing coefficient
of a simple object $V=(j^{\otimes N})\otimes\lambda$
 are given by the following formulas.

$$\langle V\rangle=\langle \lambda \rangle=
\prod_{\mathrm{cells}}\frac{[N+cn(c)]}{[hl(c)]}$$

\centerline{\begin{pspicture}[.5](0,-0.5)(1,1)
\pscurve{-}(.1,0)(.15,.45)(.45,.75)(.7,.5)(.4,.2)(.20,.38)
\pscurve{->}(.14,.62)(.1,.9)(.1,1)
 \put(-.2,0){$V$}
\end{pspicture}
$\ =\ (a^Ns)^{Ni^2+2|\lambda|}a^{|\lambda|^2}s^{N|\lambda|+2\sum_{\mathrm{cells}}cn(c)}\ $
\begin{pspicture}[.5](0,-.5)(1,1)
\psline{->}(.1,0)(.1,1)
 \put(.2,0){$V$}
\end{pspicture}
}

We denote by $\widetilde{\V}_{N,K}(\Sigma_g)$
the TQFT vector space, associated with a genus $g$ surface
$\Sigma_g$,
 for the modular category $\mathcal{C}^{PU(N,K)}$ and by
 $\tilde{d}_{N,K}(g)$
 its rank.
\begin{thm}
The rank $\tilde{d}_{N,K}(g)$
is $$\tilde{d}_{N,K}(g)=\frac{d_{N,K}(g)}{{N'}^g}$$  
\end{thm}
\begin{proof}
By Turaev's formula we have the following.
$$\tilde d_{N,K}(g)=
\left(\sum_{V\in\widetilde\Gamma_{N,K}}\langle V\rangle^2\right)
^{g-1}\sum\limits_{V\in\widetilde\Gamma_{N,K}}
\langle V\rangle^{2-2g}$$
Here  $\widetilde\Gamma_{N,K}\subset \dot{\Gamma}_{N,K}$ is a representative set of  the orbits in  $\dot\Gamma_{N,K}$ under the  order $N/\alpha$ free cyclic  action.
 Note that
this action preserves the dimension. We get
$$\tilde d_{N,K}(g)=
\left(\frac{1}{\alpha N'}\sum_{V\in\dot\Gamma_{N,K}}\langle V\rangle^2\right)
^{g-1}\frac{1}{\alpha N'}\sum\limits_{V\in\dot\Gamma_{N,K}}
\langle V\rangle^{2-2g}$$
Write $V=(1^N)^{\otimes i}\otimes \lambda$, $0\leq i<\alpha$ and 
$\lambda\in \Gamma_{N,K}$.
$$\tilde d_{N,K}(g)=
\left(\frac{1}{N'}\sum_{V\in\Gamma_{N,K}}\langle V\rangle^2\right)
^{g-1}\frac{1}{N'}\sum\limits_{V\in\Gamma_{N,K}}
\langle V\rangle^{2-2g}= \frac{d_{N,K}(g)}{{N'}^g}\ .$$
\end{proof}

The following is an integral version
of a reciprocity formula in \cite{rec}.
\begin{thm}[Level-rank duality]
One has
$\tilde{d}_{N,K}(g)=\tilde{d}_{K,N}(g)$.
\end{thm}
\begin{proof}
In the construction arising from Homfly skein theory, the parameters $N$ and $K$ play the same role, so that we can interchange rows and columns in the description of  isomorphisms classes in the modular category $\mathcal{C}^{PU(N,K)}$.
 We will get the same combinatorics as for the modular category  
$\mathcal{C}^{PU(K,N)}$.
 The result follows.
\end{proof}

\subsection{Spin decomposition of the Verlinde formula for $\mathcal{C}^{PU(N,K)}$}
Here we consider the modular category 
$\mathcal{C}^{PU(N,K)}$ in the spin case.
 This means that $d=gcd(N,K)$ is even, and that
 $N'=\frac{N}{d}$ and $K'=\frac{K}{d}$ are both odd.
 We fix the framing parameter $a$ as we did above.
\begin{thm}
Under the above hypothesis,  the category $\mathcal{C}^{PU(N,K)}$
equipped with $\varrho=(K)\otimes (1^N)$
is a  modulo $d$ spin modular category.
\end{thm}
\begin{proof}
In the modular category $\mathcal{C}^{PU(N,K)}$
the object $1^N$ and $(K)$ are  invertible
with respective orders the coprime integers $\alpha$ and $\beta$.  
 It follows that $\varrho$ is invertible with order $\alpha\beta=d$.
 The figure \ref{rho} below shows that the twist coefficient
for $\varrho$ is the product of the two twist coefficients
for $1^N$ and $(K)$ and a {\em braiding coefficient} between
 $1^N$ and $(K)$.  Using 
\cite[Prop. 1.11]{Bhec} we see that the $3$ coefficients are
respectively $(a^Ns)^N=(-1)^\beta$, $(A^Ks^{-1})^K=(-1)^\alpha$, $(a^Ns)^{2NK}=1$. The product is $-1$.

\begin{figure}
\centerline{
\begin{picture}(30,60)
\put(3,5){\includegraphics[height=40mm]{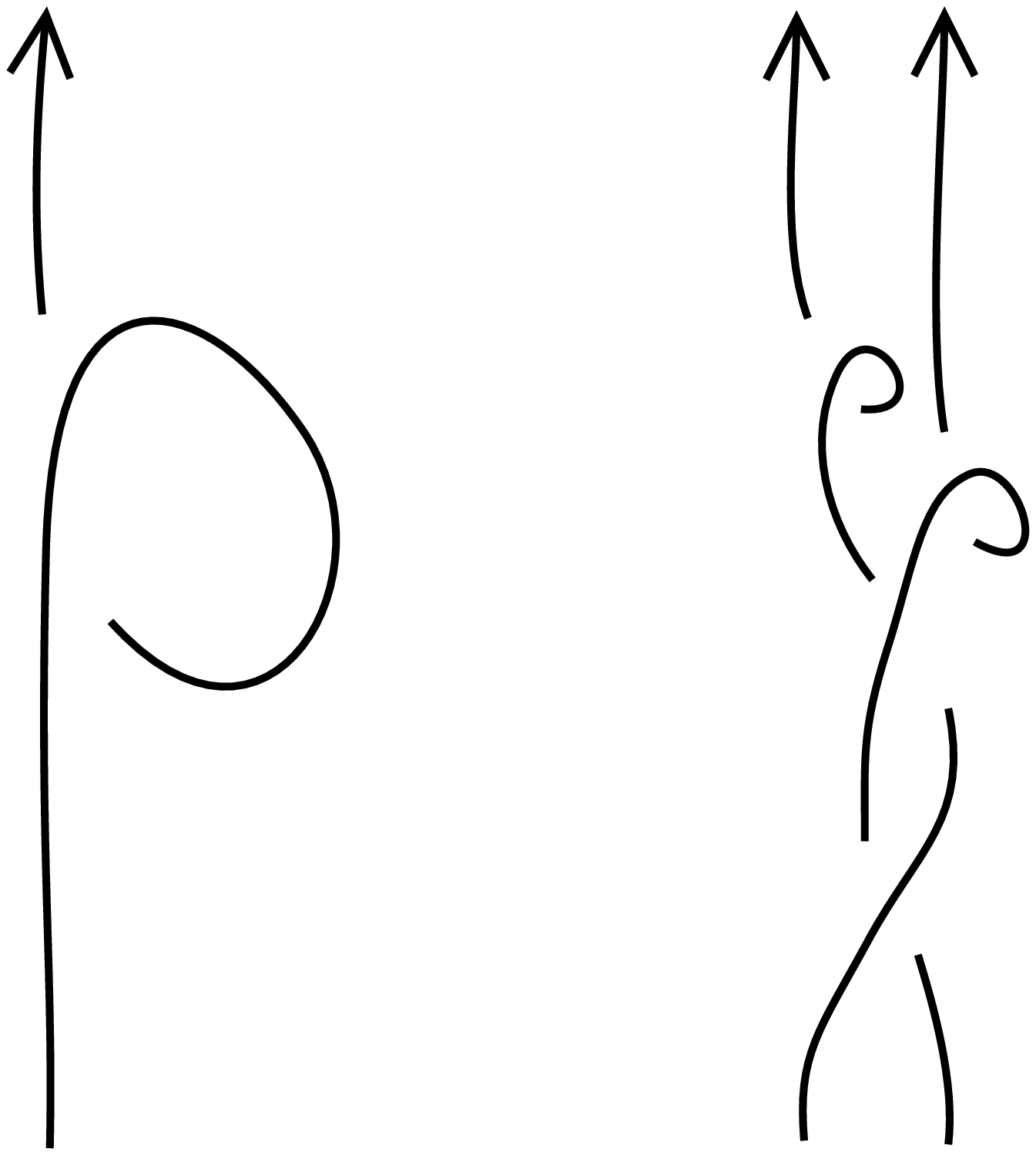}}
\put(22,25){$=$}
\put(-8,5){$V\otimes W$}
\put(27,5){$V$}
\put(38,5){$W$}
\end{picture}
}
\caption{\label{rho}Framing coefficient for a tensor product}
\end{figure}

\end{proof}

If $\sigma$ is a modulo $d$ spin structure  on the genus $g$ oriented surface
$\Sigma_g$, we denote by $\widetilde{\V}(\Sigma_g,\sigma)$
the corresponding summand and $\tilde d(g,\sigma)$ its dimension.
By applying \ref{spin}, we get.
\begin{thm}\label{spin_Atilde}
a) There exists a splitting of the Verlinde formula
$$\tilde d_{N,K}(g)=\sum\limits_{\sigma\in Spin(\Sigma_g,\mathbb{Z}/d)}
\tilde d_{N,K}(g,\sigma)\ .$$
b) The refined Verlinde formula is the following
\begin{eqnarray*} \tilde d_{N,K}(g,\sigma)&=&
\left((N+K)^{N-1}d\right)^{g-1}
  \sum\limits_{V=(1^N)^{\iota}\otimes\lambda\in\tilde{\Gamma}_{N,K}}\ \prod\limits_{\nu=1}^g
\frac{\epsilon_V(a_\nu(\sigma),b_\nu(\sigma))}
{(\sharp\, \mathrm{orb}(V))^2}\\
& &\times\prod\limits_{1\leq i<j\leq N}
  \left({2\sin{(\lambda_i-i-\lambda_j+j)\frac{\pi}{N+K}}}
\right)^{2-2g} 
\ .\end{eqnarray*}
\end{thm}
Here $\epsilon_V$ and $\sharp\mathrm{orb}(V)$
are defined as before,
 $\tilde{\Gamma}_{N,K}$ is a representative set of the quotient
of $$\dot{\Gamma}_{N,K}=
\{(j^{\otimes N})\otimes\lambda,\ j\in\mathbb{Z}/\alpha, \ \lambda\in\Gamma_{N,K}\}$$ 
under a free action of $\mathbb{Z}/\alpha N'$.
 The formula in b) can be expressed as follows.
\begin{eqnarray*} 
\tilde d_{N,K}(g,\sigma)&=&
\left((N+K)^{N-1}d\right)^{g-1}
\frac{1}{\alpha N'}
\cr & &\times
  \sum\limits_{V=(1^N)^{\iota}\otimes\lambda\in\dot{\Gamma}_{N,K}}\ \prod\limits_{\nu=1}^g
{\epsilon_V(a_\nu(\sigma),b_\nu(\sigma))}
\left(\frac{\alpha N'}{\sharp\, \mathrm{Orb}(V)}\right)^2
\cr & &\times\prod\limits_{1\leq i<j\leq N}
  \left({2\sin{(\lambda_i-i-\lambda_j+j)\frac{\pi}{N+K}}}
\right)^{2-2g} 
\ .\end{eqnarray*}

We consider now  the orbit $\mathrm{Orb}(V)$ under the action of the group
$\mathbb{Z}/\alpha\times\mathbb{Z}/N$
on $\dot{\Gamma}_{N,K}$,
where $(1,0)$ acts by
$$(1^N)^{\otimes\iota}\otimes\lambda\mapsto (1^N)^{\otimes(\iota+1)}\otimes\lambda, $$
and $(0,1)$ acts by
$$(1^N) ^{\otimes \iota}\otimes\lambda\mapsto
(1^N) ^{\otimes(\iota+\lambda_{n-1})}\otimes
((K,\lambda)-\lambda_{N-1}^N)\ .$$

\esp\noindent If $\sharp\mathrm{Orb}(V)/\alpha N'$ is even, then
$$\epsilon_V(\mathfrak{a},\mathfrak{b})=\left\{\begin{array}{l}
1\ \text{ if  $\mathfrak{a}$ and $\mathfrak{b}$ are zero modulo $|Stab(V)|$,}\\
0 \ \text{ else.}\end{array}\right.$$
If $\sharp\mathrm{Orb}(V)/\alpha N'$ is odd, then
$$\epsilon_V(\mathfrak{a},\mathfrak{b})=\left\{\begin{array}{l}
\frac{1}{2}(-1)^{\frac{2\mathfrak{a}}{|Stab(V)|}\frac{2\mathfrak{b}}{|Stab(V)|}}
\ \text{ if  $\mathfrak{a}$ and $\mathfrak{b}$ 
are  zero modulo $\frac{|Stab(V)|}{2}$,}\\
0\ \text{ else.}\\
\end{array}\right.$$
 
\section{Cohomological decomposition}\label{coho_}
In this section we will establish the decomposition in the cohomological case.

Let $d$ be an integer, and $(\mathcal{C},\varrho)$ be a modulo $d$ cohomological modular category. This means that the object $\varrho$ has order $d$
and  twist coefficient $\theta_\varrho=1$.
We deduce that the quantum dimension of $\varrho$ is $\mathfrak{\mathfrak{d}}=\pm 1$,
and\\[10pt]
\begin{equation}
\begin{pspicture}[.5](0,-0.5)(1,1)
\psline{->}(1,0)(0,1)
\psline{-}(0,0)(.35,.35)
\psline{->}(.65,.65)(1,1)
 \put(0,.3){$\varrho$} \put(.8,.3){$\varrho$}
\end{pspicture}
\ =\ 
\mathfrak{d}\ \ 
\begin{pspicture}[.5](0,-0.5)(1,1)
\pscurve{->}(0,0)(.3,.5)(0,1)
\pscurve{->}(1,0)(.7,.5)(1,1)
 \put(0,.4){$\varrho$} \put(.8,.4){$\varrho$}
\end{pspicture}
\end{equation}

\begin{equation}
\begin{pspicture}[.5](0,-0.5)(1,1)
\pscurve{->}(0,0)(.3,.5)(0,1)
\pscurve{<-}(1,0)(.7,.5)(1,1)
 \put(0,.4){$\varrho$} \put(.8,.4){$\varrho$}
\end{pspicture}
\ =\ \mathfrak{d}\ \ \ 
\begin{pspicture}[.5](0,-0.5)(1,1)
\pscurve{->}(0,0)(.5,.3)(1,0)
\pscurve{<-}(0,1)(.5,.7)(1,1)
 \put(.3,0){$\varrho$} \put(.4,.9){$\varrho$}
\end{pspicture}
\end{equation}
 
After fixing a $d$-th root of unity $\zeta$, the category
is $\Z/d$ graded.
The Kirby color  decomposes according to
this grading.
$$\Omega=\sum_{\lambda\in \Gamma}
\la \lambda \ra \lambda\ =\sum_{j\in \Z/d}\ \Omega_j$$

Using this grading we obtain the theorem below \cite{Bhec,LT}.
\begin{thm}Let $\mathcal{C}$ be a modulo $d$ cohomological modular category,
and 
\mbox{$\Omega=\sum_{j\in \Z/d}\ \Omega_j$} be the graded decomposition of the Kirby
element.
Provided \mbox{$c=(c_1,\dots,c_m)\in \Z^m$} is in the kernel of $B_L\otimes \Z/d$
the formula
$$\tau_{\mathcal{C}}^{\mathrm{coho}}(M,\sigma)=
\frac{{\langle L(\Omega_{c_1},
\dots,\Omega_{c_m})\rangle}}{\langle U_1(\Omega)\rangle^{b_+}
\langle U_{-1}(\Omega)\rangle^{b_-}}
$$
is an invariant of the surgered manifold $M={\bf S}^3(L)$
equipped with the modulo $d$ cohomology class $\sigma$ corresponding to
$c$.
\\
Moreover,
$$\forall M\ \ \tau_{\mathcal{C}}(M)=\sum_{\sigma\in Spin(M;\mathbb{Z}/d)}
\tau_{\mathcal{C}}^{\mathrm{coho}}(M,\sigma)$$
\end{thm}
Following section \ref{spinref} we get the proposition below.
 Note that here the action given by a trivial curve $\gamma$
colored with $\mathfrak{d}\phi_\gamma$ is trivial.
\begin{pro}\label{act_co}
There exists a well defined action
of the group $H_1(\Sigma,\mathbb{Z}/d)$ on $\V_\mathcal{C}(\Sigma)$, 
which maps $x=[\gamma]$ to the operator
$\psi_x=(\mathfrak{d})^{\sharp\gamma}\phi_\gamma$.
\end{pro}
Using this action we get the decomposition theorem below.
\begin{thm}\label{coho}Let $(\mathcal{C},\varrho)$ be a 
modulo $d$ cohomological modular category.\\
a) There exists a splitting of the Verlinde formula
$$ dim(\V_\mathcal{C}(\Sigma_g))=
\sum\limits_{\sigma\in H^1(\Sigma_g,\mathbb{Z}/d)}
dim(\V_\mathcal{C}(\Sigma_g,\sigma))\ .$$
b) The refined Verlinde formula is the following
\begin{eqnarray*} dim(\V_\mathcal{C}(\Sigma_g,\sigma))&=&
\langle \Omega\rangle^{g-1}
  \sum\limits_{\lambda\in\Gamma}\ 
  \la\lambda\ra^{2-2g} \times \prod\limits_{\nu=1}^g
\frac{\epsilon_\lambda(a_\nu(\sigma),b_\nu(\sigma))}
{(\sharp\, \mathrm{orb}(\lambda))^2}
\ .\end{eqnarray*}
\end{thm}
Here $(a(\sigma),b(\sigma))\in (\mathbb{Z}/d)^g\times (\mathbb{Z}/d)^g$
is given by the values of $\sigma$ on a sympleptic basis, and
$$\epsilon_\lambda(\mathfrak{a},\mathfrak{b})=\left\{\begin{array}{l}
1\ \text{ if  $\mathfrak{a}$ and $\mathfrak{b}$ are zero modulo $|Stab(\lambda)|$,}\\
0 \ \text{ else.}\end{array}\right.$$
\begin{proof}
The decomposition a) follows from the action given in Proposition \ref{act_co}.

For $\sigma\in  H^1(\Sigma_g,\mathbb{Z}/d)$, we have
$$
dim(\V_\mathcal{C}(\Sigma_g,\sigma))
=\langle \Omega\rangle^{-1-g}
\sum_\lambda \langle \lambda \rangle
B_\lambda\ ,$$
where $B_\lambda$ is the invariant of the colored link
in figure \ref{cborr}.
 Here $(a_1,b_1),\dots,(a_g,b_g)$ are given by the values
of $\sigma$ on the corresponding curves
if $\mathfrak{d}=1$, and this values plus $\frac{d}{2}$
if $\mathfrak{d}=-1$ (in this case $d$ has to be even).

The computation is done as section 3. We have
\begin{equation}\label{blambdac}B_\lambda=\langle \lambda\rangle
\prod\limits_{\nu=1}^g
\sum_{j,j'\in Stab(\lambda)} \frac{\zeta^{ja_\nu+j'b_\nu}}
{\langle \lambda\rangle^2}
\frac{\langle \Omega\rangle^2}{d^2}
\end{equation}
The formula follows.
\end{proof}
Let $d=gcd(N,K)$. If $d$ is odd, or if $d$ is even but
$\frac{NK}{d^2}$ is even, then the category
${\mathcal{C}}^{PU(N,K)}$ is a modulo $d$ cohomological modular category.

If $N=jl$, and $(-a^{N+K})^{Kj^2}=1$, then the category 
$\mathcal{C}^{SU(N,K)}$ equipped with the invertible object $\varrho=(K)^{\otimes j}$ is a 
modulo $l$ cohomological modular category.
In particular, if $N$ divides $K$, and
$N$ is odd or $\frac{K}{N}$ is even, then the category
$\mathcal{C}^{SU(N,K)}$ is a modulo $N$ cohomological modular category.

\section{Some computations}
\label{comp}
We give below some computations obtained with MuPAD \cite{Mupad}.
 Our program implements the Verlinde formulas for the categories
$\mathcal{C}^{SU(N,K)}$. The cardinality of the alcove increases rapidely,
and we obtain results only for small values of $N$,$K$.
 The function Verlinde$(N,K,g)$ gives $d_{N,K}(g)$,
and Spin\_Verl$(N,K,[\dots])$ computes $d_{N,K}(g,\sigma)$,
where  the value of $q_\sigma$ on the standard basis
is the list $[\dots]$. We know \cite{BM} that $d_{2,2}(g,\sigma)$
is $0$ or $1$ according to the Arf invariant of the spin structure.
 It would be interesting to understand the combinatorics
of the formula $d_{N,K}(g,\sigma)$ in the general case.

\begin{verbatim}
Verlinde(2,2,1);
					3
Verlinde(2,2,2);
					10
Spin_Verl(2,2,[[0,0]]);
					1
Spin_Verl(2,2,[[1,1]]);
					0
Spin_Verl(2,2,[[0,0],[1,1]]);
					0
Spin_Verl(2,2,[[1,1],[1,1]]);
					1
Verlinde(2,6,1);
					7
Verlinde(2,6,2);
					84
Spin_Verl(2,6,[[0,0]]);
					2
Spin_Verl(2,6,[[1,1]]);
					1
Spin_Verl(2,6,[[0,0],[1,1]]);
					4
Spin_Verl(2,6,[[0,0],[0,0]]);
					6
Verlinde(4,4,1);
					35
Verlinde(4,4,2);
					4680
Spin_Verl(4,4,[[0,0]]);
					3
Spin_Verl(4,4,[[1,0]]);
					2
Spin_Verl(4,4,[[1,1]]);
					2
Spin_Verl(4,4,[[2,2]]);
					2
Spin_Verl(4,4,[[0,0],[0,0]]);
					24
Spin_Verl(4,4,[[1,0],[0,0]]);
					18
Spin_Verl(4,4,[[1,0],[1,0]]);
					18
Spin_Verl(4,4,[[2,2],[0,0]]);
					20
Verlinde(6,6,1);
					462
Verlinde(6,6,2);
					30660988
Spin_Verl(6,6,[[0,0]]);
					14
Spin_Verl(6,6,[[1,0]]);
					13
Spin_Verl(6,6,[[2,0]]);
					13
Spin_Verl(6,6,[[1,1]]);
					12
Spin_Verl(6,6,[[2,2]]);
					13
Spin_Verl(6,6,[[3,0]]);
					14
Spin_Verl(6,6,[[3,3]]);
					12
Spin_Verl(6,6,[[0,0],[0,0]]);
					23718
Spin_Verl(6,6,[[1,0],[0,0]]);
					23678
Spin_Verl(6,6,[[1,0],[0,0]]);
					23624
Spin_Verl(6,6,[[2,0],[0,0]]);
					23678
Spin_Verl(6,6,[[2,2],[0,0]]);
					23678
Spin_Verl(6,6,[[3,0],[0,0]]);
					23718
Spin_Verl(6,6,[[3,3],[0,0]]);
					23648
\end{verbatim}

\bibliographystyle{amsplain}

\end{document}